\documentclass[reqno,11pt]{amsart}
\usepackage[utf8]{inputenc}

\usepackage{amsmath,amsxtra,amsbsy,enumerate, latexsym, amsfonts, amssymb, amsthm, amscd, stmaryrd}
\usepackage{hyperref}
\usepackage{graphics,epsf,psfrag}
\usepackage{bbm, dsfont}
\usepackage{mathtools}
\usepackage{xcolor}
\usepackage{float}
\usepackage{tikz}
\usetikzlibrary{patterns}
\usepackage[headings]{fullpage}
\usepackage{colonequals} 
\usepackage{microtype}
\numberwithin{equation}{section}
\allowdisplaybreaks

\definecolor{darkred}{rgb}{0.7,0.1,0.1}

\definecolor{darkgreen}{rgb}{0.1,0.7,0.1}

\newcommand{\bbE}{{\ensuremath{\mathbb E}} }

\newcommand{\bbN}{{\ensuremath{\mathbb N}} }

\newcommand{\bbP}{{\ensuremath{\mathbb P}} }

\newcommand{\bbR}{{\ensuremath{\mathbb R}} }

\newcommand{\bbZ}{{\ensuremath{\mathbb Z}} }

\renewcommand{\epsilon}{\varepsilon}

\newcommand{\ga}{\alpha}
\newcommand{\gb}{\beta}

\newcommand{\gep}{\varepsilon}       

\newcommand{\gG}{\Gamma}

\newcommand{\go}{\omega}

\newcommand{\gO}{\Omega}
\newcommand{\gl}{\lambda}
\newcommand{\gL}{\Lambda}
\newcommand{\gs}{\sigma}

\newcommand{\cG}{{\ensuremath{\mathcal G}} }

\newcommand{\cA}{{\ensuremath{\mathcal A}} }
\newcommand{\cB}{{\ensuremath{\mathcal B}} }
\newcommand{\cF}{{\ensuremath{\mathcal F}} }

\newcommand{\cC}{{\ensuremath{\mathcal C}} }
\newcommand{\cN}{{\ensuremath{\mathcal N}} }
\newcommand{\cL}{{\ensuremath{\mathcal L}} }
\newcommand{\cT}{{\ensuremath{\mathcal T}} }

\newcommand{\bP}{{\ensuremath{\mathbf P}} }

\newcommand{\bE}{{\ensuremath{\mathbf E}} }

\DeclareMathOperator*{\infess}{ess\, inf}

\newcommand{\ind}{\mathbf{1}}

\newcommand{\lint}{\llbracket}
\newcommand{\rint}{\rrbracket}

\newtheorem{theorem}{Theorem}[section]
\newtheorem{lemma}[theorem]{Lemma}
\newtheorem{proposition}[theorem]{Proposition}

\newcommand{\RN}[1]{%
  \textup{\uppercase\expandafter{\romannumeral#1}}%
}

\newcommand{\Var}{\mathrm{Var}}
\newcommand{\TV}{\mathrm{TV}}
\newcommand{\Gap}{\mathrm{gap}}

\newcommand{\dd}{\mathrm{d}}
\newcommand{\Mix}{\mathrm{mix}}

\renewcommand{\tilde}{\widetilde}

\newcommand{\maxtwo}[2]{\max_{\substack{#1 \\ #2}}} 
\newcommand{\suptwo}[2]{\sup_{\substack{#1 \\ #2}}} 
\newcommand{\inftwo}[2]{\inf_{\substack{#1 \\ #2}}} 
\newcommand{\sumtwo}[2]{\sum_{\substack{#1 \\ #2}}} 

\makeatletter
\def\captionfont@{\footnotesize}
\def\captionheadfont@{\scshape}

\long\def\@makecaption#1#2{%
  \vspace{2mm}
  \setbox\@tempboxa\vbox{\color@setgroup
    \advance\hsize-6pc\noindent
    \captionfont@\captionheadfont@#1\@xp\@ifnotempty\@xp
        {\@cdr#2\@nil}{.\captionfont@\upshape\enspace#2}%
    \unskip\kern-6pc\par
    \global\setbox\@ne\lastbox\color@endgroup}%
  \ifhbox\@ne 
    \setbox\@ne\hbox{\unhbox\@ne\unskip\unskip\unpenalty\unkern}%
  \fi
  \ifdim\wd\@tempboxa=\z@ 
    \setbox\@ne\hbox to\columnwidth{\hss\kern-6pc\box\@ne\hss}%
  \else 
    \setbox\@ne\vbox{\unvbox\@tempboxa\parskip\z@skip
        \noindent\unhbox\@ne\advance\hsize-6pc\par}%
\fi
  \ifnum\@tempcnta<64 
    \addvspace\abovecaptionskip
    \moveright 3pc\box\@ne
  \else 
    \moveright 3pc\box\@ne
    \nobreak
    \vskip\belowcaptionskip
  \fi
\relax
}
\makeatother

\newcommand{\gap}{\mathrm{gap}}

\def\writefig#1 #2 #3 {\rlap{\kern #1 truecm
\raise #2 truecm \hbox{#3}}}

\title[ASEP in random environment]{ Mixing time for the asymmetric simple exclusion process in a random environment}

\author[Hubert Lacoin]{Hubert Lacoin}
 \address{Hubert Lacoin \hfill\break
IMPA\\
Estrada Dona Castorina, 110\\ Rio de Janeiro 22460-320 \\ Brazil.}
\email{lacoin@impa.br}

\author[Shangjie Yang]{Shangjie Yang}
 \address{Shangjie Yang \hfill\break
IMPA\\
Estrada Dona Castorina, 110\\ Rio de Janeiro 22460-320 \\ Brazil.}
\email{yashjie@impa.br}


\begin{document}


\begin{abstract}
We consider the simple exclusion process in the integer segment $\lint 1, N \rint$ with $k\le N/2$ particles and  spatially inhomogenous jumping rates.  A particle at site $x\in \lint 1, N \rint$ jumps to site $x-1$ (if $x\ge 2$) at rate $1-\go_x$ and to  site $x+1$ (if $x \le N-1$) at rate $\go_x$ if the target site is not occupied. The sequence  $\go=(\go_x)_{ x \in \bbZ}$ is chosen by IID sampling from a probability law whose 
support is bounded away from zero and one (in other words the random environment satisfies the uniform ellipticity condition).
We further assume 
$\bbE[ \log  \rho_1 ]<0$ where $\rho_1:= (1-\go_1)/\go_1$,
which implies that our particles  have a tendency to move to the right.
We prove that the mixing time of the exclusion process in this setup grows like a power of $N$. More precisely,  for the exclusion process with $N^{\beta+o(1)}$ particles where $\beta\in [0,1)$, we have in the large $N$ asymptotic  
$$   N^{\max\left(1,\frac 1\gl, \beta+ \frac 1 {2\gl}\right)+o(1)} \le  t_{\Mix}^{N,k}  \le N^{C+o(1)}$$
where  $\lambda>0$ is such that $\bbE[\rho_1^{\gl}]=1$ ($\gl=\infty$ if the equation has no positive root) and $C$ is a constant which depends on the distribution of $\go$. 
We conjecture that our lower bound is sharp up to sub-polynomial correction.

\end{abstract}

\maketitle


\date{\small\today}
\maketitle

\section{Introduction}\label{sec:introd}

\subsection{Overview}

From the viewpoint of Probability and Statistical Mechanics,
 the simple exclusion process is  one of the simplest interacting particle system. It is a reasonable toy model to describe the relaxation of a low density gas and we refer to \cite[Chapter VIII.6]{liggett2012interacting} for a historical introduction. Its relaxation to equilibrium has been the object of extensive study under  a variety of perspective: Hydrodynamic limits  \cite{Rost81, KOV89, rez91}, Relaxation Time \cite{DSC93, Quastel92} $\log$-Sobolev inequalites \cite{Yau97} and Mixing Time \cite{benjamini2005mixing, Morris06} (the list of references is very far from exhaustive).

 \medskip

All the above mentioned works are concerned with the exclusion in an \textit{homogeneous} medium and
a small modification of this setup can lead to a drastic change of the pattern of relaxation, see for instance \cite{FGS16, FN17} (and references therein) for the phenomenology induced by the change of the jump rate on a single bond.
The \textit{disordered} setup, where the jump rate of the particles is random and varies in space fostered interest only more recently, see for instance \cite{faggionato2008random, FRS19, schmid2019mixing}. 

\medskip

In the present paper, we are interested in the case IID site disorder on a one dimensional segment, in particular in the case where the local drift felt by particles has a non-constant sign. For the system to reach equilibrium, individual particles need to travel on macroscopic distances and in particular  have to fight against drift in some regions. This phenomenon, also present in the case of the random walk in a random environment \cite{gantert2012cutoff, kesten1975limit}, induces a slower mixing than in the constant nonzero bias case, as was proved in \cite{schmid2019mixing}.
Our objective is to quantify further this slow down of the mixing time.

\medskip

In order to estimate the mixing time of the  disordered exclusion process, we need to understand in details how these regions with unfavorable drift -- which we refer to as \textit{traps} -- affect the pattern of relaxation to equilibrium. We make two important steps toward this objective: 
\begin{itemize}
 \item We prove that the mixing time grows at most like a power of $N$ (the upper bound we prove displays a non optimal exponent). 
 \item We obtain a lower bound on the mixing time, which we conjecture to be optimal, and which allows to identify, depending on the parameters of the system, which is the main factor that slows down the mixing. 
\end{itemize}
More precisely, our proof of the lower bound shows that the mixing time can be bounded from below by three different mechanisms:
\begin{itemize}
 \item [$(i)$] Particles cannot move faster than ballistically, so that the mixing time is  at least of order $N$ which is the length of the system.
 \item [$(ii)$] The particles may remain trapped in potential wells which are created by the environment (see the definition \ref{defvpot}), so that the mixing time is at least of order $e^{\Delta V}$ where $\Delta V$ is the height of the worse potential well in the system.
 \item[$(iii)$] The potential wells also limit the flow of particles through the system which is at most of order $e^{-\Delta V/2}$. For this last reason, the mixing time is at least of order $k e^{\Delta V/2}$ when $k$ is the number of particles in the system.
\end{itemize}
While the two first limitations $(i)$ and $(ii)$ follow from early studies of one dimensional random walk in a random environment \cite{kesten1975limit} and have already been used to determine its mixing time  \cite{gantert2012cutoff}. The third limitation is specific to systems with many particles, and to our knowledge, had not been identified so far. It creates a third phase in the conjectured mixing time diagram (see Figure \ref{fig:dynamicdiagram}).

\subsection{The exclusion process in a random environment}

Let us introduce formally the random process whose study is the object of this paper. 
The exclusion process on the segment $\lint 1,N\rint$  with $k$ particles and $1\le k \le N/2$ is a Markov process that can informally be described as follows.

 \begin{itemize}
 \item [(A)]
 Each site is occupied by at most one particle (we refer to this constraint as \textit{the exclusion rule)}. Therefore at all time there are $k$ occupied sites and $N-k$ empty sites.
 \item [(B)] Each of the $k$ particles performs a random walk on the segment, independently of the others, except that any jump that violates the exclusion rule is  cancelled.
    \end{itemize}
More precisely, we want to consider the case  exclusion process in a \textit{random environement} where the jump rates of the particles are specified by sampling an IID 
sequence of random variables $\go=(\go_x)_{x\in \bbZ}$, and the transition rates are given by 
\begin{equation}\begin{cases}\label{qngo}
                  q^{\go}_N(x,x+1)= \go_{x}\ind_{\{ x\le  N-1\}},&\\
                  q^{\go}_N(x,x-1)= (1-\go_{x})\ind_{\{ x\ge 2\}},&\\
                  q^{\go}_N(x,y)=0   & \text{ if } y \notin \{x-1,x+1\}.
                \end{cases}
\end{equation}
The random walk with transitions $q^{\go}_N$ which corresponds to the case $k=1$  is an extensively studied process, usually 
referred to as Random Walk in a Random Environment (RWRE).
The RWRE on the full line $\bbZ$ was first studied by Solomon in \cite{solomon1975random} who established a criterion for recurrence/transience. 
The limit law of the random walk in a random environment is studied by  Kesten \textit{et al.} in \cite{kesten1975limit}  when the random walk is transient, and by Sinai in \cite{sinai1982limit} when the random walk is recurrent (we refer to \cite{sznitman2004topics, zeitouni2004part}  for complete introductions to this research field).

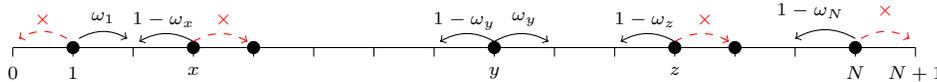
\begin{figure}[h]
 \centering
   \begin{tikzpicture}[scale=.4,font=\tiny]
     \draw (23,-1) -- (53,-1);
     
     \foreach \x in {23,25, 27,...,51,53} {\draw (\x,-1.3) -- (\x,-1);}
     \draw[fill] (51,-1) circle [radius=0.2];   
     \node[below] at (25,-1.3) {$1$};
     \node[below] at (51,-1.3) {$N$};
     \node[below] at (23,-1.3) {$0$};
     \node[below] at (53,-1.3) {$N+1$};

       \draw  (24.8, -0.7)  edge[bend right,->, dashed,red](23.2, -0.7);
       \node[below, red] at (24, 0.5) {$\times$};
       \draw  (25.2, -0.7)  edge[bend left,->](26.8, -0.7);
       \node[below] at (26, 0.5) {$\go_1$};

       \draw[fill] (25,-1) circle [radius=0.2];  
       \node[below] at (29, -1.3) {$x$};
       \node[below] at (28, 0.5) {$1-\go_x$};
       \node[below,red] at (30, 0.5) {$\times$};    

       \draw (51.2, -0.7) edge[bend left,->, dashed,red] (52.8, -0.7);
       \draw (51, -0.7) edge[bend right,->] (49, -0.7);
       \node[below] at (49.5,0.8) {$1-\go_{N}$};
       \node[below,red] at (52,.8) {$\times$};

       \draw[fill] (29,-1) circle [radius=0.2];
   
       \draw (29, -0.8) edge[bend right,->] (27.2, -0.8);
       \draw (29, -0.8) edge[bend left,dashed,->, red] (30.8, -0.8);
       \draw[fill] (31,-1) circle [radius=0.2]; 
       
       \draw (39, -0.8) edge[bend left,->] (40.8, -0.8);   
       \draw (39, -0.8) edge[bend right,->] (37.2, -0.8);      
       \node[below] at (39,-1.3) {$y$};
       \node[below] at (40,0.5) {$\go_y$};
       \node[below] at (38,0.5) {$1-\go_y$};

       \draw[fill] (39,-1) circle [radius=0.2];
         
       \draw (45, -0.8) edge[bend right, ->] (43.2, -0.8);  
       \draw (45, -0.8) edge[bend left,->,dashed,red] (46.8, -0.8);      
       \node[below, red] at (46,0.5) {$\times$};
       \node[below] at (44,0.5) {$1-\go_z$};
       \node[below] at (45,-1.3) {$z$};
   
       \draw[fill] (45,-1) circle [radius=0.2];   
       \draw[fill] (47,-1) circle [radius=0.2];       
 \end{tikzpicture}
   \caption{ A graphical representation of the simple exclusion process in the segment $\lint 1, N \rint$ and environment $\go=(\go_x)_{x \in \bbZ}$: a bold circle represents a particle, and the number above every arrow represents the jump rate while a red $"\times"$ represents an inadmissible jump.} \label{fig:asep} 
 \end{figure}

\noindent    

 We are interested in the following quantitative question: How long does the system need to relax to equilibrium, forgetting the information of its initial configuration in the sense of total-variation distance?  More precisely we are interested in the asymptotic in the limit when $k,N\to\infty$ of this \textit{total-variation mixing time}. 
 This question has been extensively studied in the case where the sequence $\go=(\go_x)_{x \in \bbZ}$ is constant, which we refer to as the \textit{homogeneous environment} case:
\begin{itemize}
\item[(1)]
When $\go_x \equiv \frac{1}{2}$, Wilson in \cite{wilson2004mixing} showed that the system takes time of order $N^2 \log \min(k, N-k)$ and later Lacoin in \cite{lacoin2016mixing} proved that the lower bound in \cite{wilson2004mixing} is sharp. 
 
 \item[(2)]
 When $\go_x \equiv p \neq \frac{1}{2}$, Benjamini \textit{et al.} in \cite{benjamini2005mixing} told that the system takes time of order $N$, and later Labb\'{e} and Lacoin in \cite{labbe2016cutoff} provided the exact constant. 
 
\item[(3)] 
 The case  $\go_x \equiv p_N=\frac{1}{2}+\gep_N$ with $\lim_{N\to \infty}\gep_N=0$ is studied by Levins and Peres in \cite{levin2016mixing},  Labb\'{e} and Lacoin in \cite{labbe2018mixing}.
\end{itemize}

\medskip

 From the results mentioned above, for homogeneous environment the system takes time  at least of order $N$ and  most of order $N^2 \log N$ to relax to equilibrium. 
However,
when the sequence $\go=(\go_x)_{x \in \bbZ}$ is chosen by independently sampling a nondegenerate common law, the system can  exhibit a very different behavior because the random environment 
can create  wells of potential which trap particles (see Equation \eqref{defvpot} below for a definition of the potential associated to $\go$).

\medskip

  Gantert and Kochler has studied the mixing time problem when $k=1$ (and transient environment)  in \cite{gantert2012cutoff} for random environment and  identified the mixing time, which is related with the depth of the deepest trap  and may be much larger than $N^2 \log N$. Schmid \cite{schmid2019mixing} studied the question in the case of a  positive density of particles, when the environment is ballistic to the right, (that is, when the random walk  is transient with positive speed) and provided bounds for the mixing time, showing in particular that  the mixing time is of  order larger than $N$ as soon as 
  the local drift (which is equal to $2\go_x-1$) is not uniformly bounded from below by a positive constant and is larger than $N^{1+\delta}$ for some $\delta>0$ when some sites can display negative drift
($\bbP[\go_x<1/2]>0$).

\medskip

 In our study we focus on the case of random environments  which are  such that the random walk is transient (the case of recurrent environment is quite different and should be considered separately). In that setup, the results in \cite{schmid2019mixing}
leave several questions open, among which the following ones:
\begin{itemize}
\item [(A)] Is the mixing time always bounded from above by a power of $N$?
\item [(B)] If this is the case, for the exclusion process with $k_N=N^{\beta}$ particles and $\beta\in (0,1)$,
can one identify an exponent $\nu>0$ (depending on $\beta$ and the distribution) which is such that the mixing time is of order $N^{\nu}$?
\end{itemize}
We provide a positive answer to question $(A)$ by proving an upper bound on the mixing time which grows like a  power of $N$. This upper bound is achieved by using a censoring procedure which allows to transport particles one by one to their equilibrium positions.
Concerning question $(B)$, we provide a new lower bound on the mixing time which we believe to be optimal and provide a conjecture concerning the value of $\nu$. The bound is based on an analysis of the effect of the deepest trap on the particle flow through the system. Significant technical obstacles prevented us from obtaining a matching upper bound.

\section{Model and result}\label{sec:modres}

\subsection{An introduction to  Random Walk in a Random  Environment $\go$}

Let us recall the definition for random walk in a random environment. Given
 $\omega=(\omega_x)_{x\in \bbZ}$ a sequence with values in $(0,1)$, the random walk in the environment $\go$ is the continuous time Markov chain  on $\bbZ$ whose transition rates are given by
 \begin{equation}\label{defqgo}
  \begin{cases}
   q^{\go}(x,x+1)= \go_x,\\
   q^{\go}(x,x-1)= 1-\go_{x},&\\
   q^{\go}(x,y)=0 &\quad \text{ if } |x-y|\ne 1.
  \end{cases}
 \end{equation}
We let $(X_t)_{t \ge 0}$ denote the random walk in environment $\go$ and initial condition $0$ (we let $Q^{\go}$ denote the corresponding law).
This process has been extensively studied in the case where 
 $\go=(\go_x)_{x \in \bbZ}$ is (the fixed realization of) a sequence of IID random variables (we will use $\bbP$ and $\bbE$ denote  the  associated law and expectation respectively), and we refer to \cite{sznitman2004topics,zeitouni2004part} for recent reviews.
 
 \medskip
 
 Simple criteria have been derived on the distribution of $\go$ as necessary and/or sufficient conditions for recurrence/transience, ballisticity etc... Even though most of the results are valid in a more general setup, for the sake of simplicity let us assume in the discussion that 
 the variables $(\go_x)_{x \in \bbZ}$ are bounded away from $0$ and $1$, that is, for some $\alpha\in (0,1/2)$ we have 
\begin{equation}\label{unifeliptic}
 \bbP( \go_1\in [\alpha,1-\alpha])=1.
\end{equation}
Setting $\rho_x \colonequals (1-\omega_x)/ \omega_x$, 
it has been proved in \cite{solomon1975random} that 

\begin{equation}
 \begin{cases}
  \bbE[\log \rho_1]=0 \Rightarrow X_t \text{ is reccurent under $Q^{\go}$ , $\bbP$-a.s.},\\
 \bbE[\log \rho_1]\ne 0  \Rightarrow X_t \text{ is transient under $Q^{\go}$ , $\bbP$-a.s.}
 \end{cases}
\end{equation}
More precisely in the second case we have with probability one $\lim_{t\to \infty} X_t=\infty$ (resp. $-\infty$) if  $\bbE[\log \rho_1]<0$   (resp. $\bbE[\log \rho_1]>0$).

\medskip

\noindent When transience holds, the rate at which $X_t$ goes to infinity has also been identified in \cite{kesten1975limit}. It can be expressed in terms of a simple parameter of the distribution $\go$, yielding in particular a necessary and sufficient condition for ballisticity.
Let us assume that $\bbE[\log \rho_1]<0$, and 
set 
$$\gl= \gl_{\bbP}:=\inf \{ s>0, \bbE[\rho^{s}_1]\ge 1 \} \in (0,\infty]. $$ 
It has been proved in \cite{kesten1975limit} that if $\gl>1$ then there exists $\vartheta_{\bbP}>0$ such that  
\begin{equation}
 \lim_{t\to \infty} \frac{X_t}{t}=\vartheta_{}
\end{equation}
and that if $\gl\in (0,1]$ then 
\begin{equation}
 \lim_{t\to \infty} \frac{\log (X_t)}{ \log t}=\gl.
\end{equation}

\subsection{The Simple Exclusion process in an environment $\go$}

\subsubsection{Definition}
Given a sequence 
 $\omega=(\omega_x)_{x\in \bbZ}$ taking values in $(0,1)$,
$N\ge 2$  and $1\le k \le N-1$,
the simple exclusion process in a random environment on the line segment $\lint 1, N \rint$ (we use the notation $\lint a, b \rint \colonequals [a, b] \cap \bbZ$) with $k$ particles is a Markov process on the space
\begin{equation}\label{statespace}
\Omega_{N,k} \colonequals \left\{\xi \in \left\{0, 1\right\}^N: \sum_{x=1}^{N}\xi(x)=k \right\}.
\end{equation}
The $1$s' are denoting particles while $0$s' correspond to empty sites. 
It can  be informally described as follows: each of the $k$ particles performs independently a random walk with transitions given by $q^{\go}$ in \eqref{defqgo}, with the constraints that particles must remain in the segment and 
each site can be occupied by at most one particle. All transitions that would make this constraint violated (that is, a particle tries to jump either on $0$, $N+1$ or an already occupied site) are cancelled.

\medskip

More formally we let $\xi^{x,y}$ is the configuration obtained by swapping the values of $\xi$ at sites $x$ and $y$ of the configuration $\xi$,   more formally defined by
\begin{equation}\label{swap2sites}
  \forall z \in \lint 1, N \rint, \quad \xi^{x,y}(z)= \xi(z)\ind_{\lint 1,N\rint\setminus \{ x,y\}}+\xi(x)\ind_{\{y\}}+\xi(y)\ind_{\{x \}} .
\end{equation}
The simple exclusion process in environment $\go$ is the Markov process with transition rates given by
\begin{equation}\label{jumprateSEP}
\begin{split}
 r^{\go}(\xi,\xi^{x,x+1})&:= \begin{cases} \omega_x   &\text{ if } \xi(x)=1 \text{ and } \xi(x+1)=0,\\
                        1-\omega_{x+1}   \quad &\text{ if } \xi(x+1)=1 \text{ and } \xi(x)=0,\\
                     \end{cases}   \quad \text{ for } x\in \lint 1,N-1\rint\\
 r^{\go}(\xi,\xi')&:=0     \quad \text{ in all other cases}.                    
\end{split}
\end{equation}
Equivalently  the generator of the process is defined for $f: \gO_{N,k}\to \bbR$ by
\begin{equation}\label{generator}
\mathcal{L}^{\go}_{N,k}(f)(\xi) \colonequals \sum_{x=1}^{N-1} r^{\go}(\xi,\xi^{x,x+1})
\big[f(\xi^{x,x+1})-f(\xi)\big].
\end{equation}
The chain is ergodic and reversible. In order to give a simple compact expression for the equilibrium measure, let us introduce the random potential $V^{\go}:   \bbN \to  \bbR$ defined as follows, $V^{\go}(1) \colonequals 0$ and for $x \ge 2$ 
\begin{equation}
 V^{\go}(x):= \sum_{y= 2}^x \log \left( \frac{1-\go_y}{\go_{y-1}}\right).
\end{equation}
With a small abuse of notation, we extend $V^{\go}$ to a function of $\gO_{N,k}$. This extension is obtained by summing the value of $V^{\go}$ among the positions of the particles in the configuration $\xi$:
\begin{equation}\label{defvpot}
 V^{\go}(\xi):= \sum_{x=1}^N V^{\go}(x)\xi(x).
\end{equation}
We consider the probability measure  $\pi_{N,k}^{\omega}$ defined by
\begin{equation}\label{invarmeas}
\pi_{N, k}^{\omega}(\xi)\colonequals \frac{1}{Z_{N, k}^{\omega}} e^{-V^{\go}(\xi)}  \quad  \text{  with }   \quad Z_{N, k}^{\omega}= \sum_{\xi \in \gO_{N,k}} e^{-V^{\go}(\xi)}.
\end{equation}
It is immediate to check by inspection that  $\pi_{N, k}^{\omega}$ satisfies the detailed balance condition for $\mathcal{L}^{\go}_{N,k}$, and thus that it is the unique invariant probability measure on $\gO_{N,k}$.

If $\xi \in \gO_{N,k}$, we let 
 $(\sigma^{\xi}_t)_{t\ge 0}$ denote the Markov chain with initial condition $\xi$.  We are going to provide a construction $(\sigma^{\xi}_t)_{t \ge 0}$ for all $\xi\in \gO_{N,k}$ on a common probability space in Section \ref{subsec:cancoupling}, and we use  $\bP$ and $\bE$ for the  corresponding  probability law and expectation respectively. We let $(P_t)_{t\ge 0}$ (the dependence in $\go$, $N$, $k$ is omitted in the notation to keep it light) denote the corresponding Markov semi-group and set $P_{t}^{\xi}:= \bP ( \sigma^{\xi}_t \in \cdot)=P_t(\xi,\cdot)$ to be the marginal distribution of  $(\sigma^{\xi}_t)_{t \ge 0}$ at time $t$. 
 \subsubsection{Mixing time and spectral gap}
  In a standard fashion, we set  the total variation-distance to equilibrium at time $t$ to be
\begin{equation}\label{def:mixdist}
d_{N, k}^{\omega}(t) \colonequals \max_{\xi \in \Omega_{N,k}} \Vert P_t^{\xi}-\pi_{N, k}^{\omega}\Vert_{\TV}
\end{equation}
where $\Vert \nu_1- \nu_2 \Vert_{\TV} \colonequals \sup_{A \subset \Omega_{N,k}} \vert \nu_1(A)-\nu_2(A)\vert$ denotes the total variation between two probability measures $\nu_1$, $\nu_2$ on  $\Omega_{N,k}$.
Since the Markov chain is irreducible, we know that (cf. \cite[Theorem 4.9]{LPWMCMT})
\begin{equation}\label{dtozero}
\lim_{t \to \infty}d_{N, k}^{\omega}(t)=0.
\end{equation}
We are interested in having quantitative statements related to the convergence \eqref{dtozero},
and for this reason we want to evaluate the mixing time and spectral gap of the chain (see \cite{LPWMCMT} for a motivated and thorough introduction to these notions).
For  $\epsilon \in (0, 1)$, let the $\gep$-mixing time of the chain be defined by
\begin{equation}\label{def:mixtime}
t_{\Mix}^{N,k,\go}(\epsilon)\colonequals \inf \Big\{t\geq 0: d_{N, k}^{\omega}(t) \leq \epsilon \Big \}.
\end{equation}
By convention, we simply write $t_{\Mix}^{N,k,\go}$ when $\gep=1/4$.
The spectral gap of the chain $\Gap_{N, k}^{\omega}$, in our context, can be defined as 
the smallest non-zero eigenvalue of $-\mathcal{L}^{\go}_{N,k}$.
It can be shown using a spectral decomposition (see for instance  \cite[Corollary 12.7]{LPWMCMT}) to determine the asymptotic rate of convergence of $d_{N, k}^{\omega}$ as
\begin{equation}
 \lim_{t\to \infty} \frac{1}{t} \log d_{N, k}^{\omega}(t)=-\Gap_{N, k}^{\omega}.
\end{equation}
The mixing time and spectral gap are related to one another by the following relation valid for $\gep\in (0,1/2)$ (cf. \cite[Theorems 12.4 and 12.5]{LPWMCMT})

\begin{equation}\label{gapmixrelation}
 \frac{1}{  \Gap_{N, k}^{\omega} }   \log \left(\frac 1 {2\gep}\right)    \le t_{\Mix}^{N,k,\go}(\epsilon) \le  \frac{1}{  \Gap_{N, k}^{\omega} } 
 \log\left( \frac{1}{\gep \pi_{\min}}\right)
\end{equation}
where 
$$\pi_{\min}= \min_{\xi\in \gO_{N,k}} \pi^{\go}_{N,k}(\xi).$$

 \subsection{Results}
 
 The main object of the paper is the study of the exclusion process in an IID environment. On the way to our main result, we also prove bounds on the mixing time which are valid for any realization of $\go$, and which we present first. 
 
 \subsubsection{Universal bounds for the mixing time on the exclusion process}

 We assume without loss of generality (by symmetry) that $k\le N/2$.
We prove that  the mixing time grows at least linearly  with the size of the system and at most exponentially. Both results are in a sense optimal (see the discussion in Section \ref{sec:blabla}. below).

\begin{proposition}\label{prop:ulb}
For any $k\in \lint 1,N/2\rint$ and $N\ge 2$, for any  $(\go_x)_{x\in \bbZ}$
 we have 
 \begin{equation}\label{mix1}
  t_{\Mix}^{N,k,\go} \ge \frac{1}{16}{N}.
 \end{equation}
Furthermore, if $k_N$ is a sequence such that 
\begin{equation}\label{totheinfinite}
k_N\le N/2 \text{ and  } \lim_{N\to \infty} k_N=\infty,
\end{equation}
we have for any $\gep> 0$, for $N\ge N_0(\gep)$ sufficiently large  for any $(\go_x)_{x\in \bbZ}$
\begin{equation}\label{mix2}
 t_{\Mix}^{N,k_N,\go}(1-\gep)\ge  \frac{1}{30}N.
\end{equation}
\end{proposition}
\noindent For the upper bound, we require an assumption similar to \eqref{unifeliptic}, that is 
\begin{equation}\label{detele}
 \forall x\in \bbZ, \quad \go_x\in [\alpha,1-\alpha].
\end{equation}

\begin{proposition} \label{prop:lwbdgap}
For any sequence $(\go_x)_{x\in \bbZ}$ satisfying \eqref{detele} all $N\geq 2$ and all $k \in \lint 1, N/2\rint$,  
we have
\begin{equation}
\Gap_{N, k}^{\omega}\geq  \ga  N^{-2} |\gO_{N,k}|^{-1}\left(\frac{1-\alpha}{\alpha}\right)^{ - N/2},
\end{equation}
and as a consequence 
 for all $\gep \in (0,1/2)$
\begin{equation}\label{dtmmix} t_{\Mix}^{N,k,\go}(\gep) \le \ga^{-1} N^{2} \vert \gO_{N,k} \vert \left( \frac{1-\ga}{\ga} \right)^{N/2}\left( \log \vert \gO_{N,k} \vert +Nk \log \frac{1-\ga}{\ga}  - \log \gep \right).
\end{equation}

\end{proposition}

\subsubsection{Mixing time for the exclusion process in a random environment}
 
Let us now introduce our main results concerning the exclusion process in a random environment.
 We assume that \eqref{unifeliptic} holds,
\begin{equation}\label{sumpkrho}
 \bbE[\log \rho_1 ]<0   \quad \text{ and }  \quad  1\le k\le N/2.
\end{equation}
Using the various symmetries of the the system (between left and right, particles and empty sites...), assumption \eqref{sumpkrho} entails almost no-loss of generality, and the only case being left aside is that of a recurrent environment (that is  $\bbE[\log \rho_1 ]=0$).
We are also going to consider that $\gl_{\bbP}<\infty$, this corresponds to saying that 
$\bbP[\  \go_1<1/2 ]>0$ (the case $\bbP[\  \go_1\ge 1/2 ]=1$ is discussed in the next section).

\medskip
In order to bet a better intuition on the result, let us provide a description of the equilibrium measure.
We introduce the event $\cA_r \subset \gO_{N,k}$  that the leftmost particle and rightmost empty site are at a distance smaller than $2r$ of their respective maximal and minimal possible values:
\begin{equation}\label{event:Ar}
 \cA_r:=\left\{ \xi \in  \gO_{N,k} \ : \ \forall x\in \lint 1, N-k-r\rint, \ \xi(x)=0\ ;  \ \forall x\ge N-k+r, \xi(x)=1 \right\}.
\end{equation}
The following result tells us that the mass of  $\pi_{N,k_N}$ is essentially concentrated at a finite distance of the configuration  $\xi_{\max}$  with all $k$ particles packed to the right  (see \eqref{minmaxconf}).

\begin{lemma}\label{lema:equilleftright}
 Under the assumptions \eqref{detele} and \eqref{sumpkrho}, 
  for all $N$ sufficiently large we have w.h.p.
 \begin{equation}
 \lim_{r\to \infty}\inftwo{N\ge 1}{k \in \lint 1, N/2\rint} \bbE\left[ \pi_{N,k_N}^{\go}\left( \cA_r\right) \right]=1.
 \end{equation}
 
\end{lemma}

Our  first main result is that  if the environment satisfies the assumptions \eqref{unifeliptic} and \eqref{sumpkrho}, the system relaxes to equilibrium  in polynomial time, or in other words that $t_{\Mix}^{N,k,\go}$ grows like a power of $N$ with an explicit upper bound on the growth exponent.
 In order to describe our explicit bound, we need to introduce the function $F$ which is the 
$\log$-Laplace transform of $\log \rho_1$ that is 
\begin{equation}\label{loglaplace}
F(u) \colonequals \log \bbE \left[ \rho_1^{u} \right].
\end{equation}
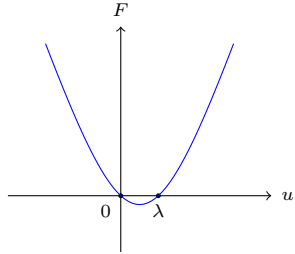
\begin{figure}[h]
\centering
  \begin{tikzpicture}[scale=.5,font=\tiny]
   \draw[->] (-3, 0) -- (4, 0) node[right] {$u$};
  \draw[->] (0, -1.5) -- (0, 4.5) node[above] {$F$};
  \node at (1, -0.4){$\gl$};
  \node at (-0.4, -0.4){$0$};

  \filldraw[black, thick] (1,0) circle (1pt);
  \filldraw[black, thick] (0,0) circle (1pt);
 
\draw[blue] (-2.,4.0464 )--(-1.9,3.78713 )--(-1.8,3.53043 )--(-1.7,3.27665 )--(-1.6,3.0262 )--(-1.5,2.77952 )--(-1.4,2.5371)--(-1.3,2.29948 )--(-1.2,2.06729 )--(-1.1,1.84117 )--(-1.,1.62186 )--(-0.9,1.41015 )--(-0.8,1.20688 )--(-0.7,1.01296 )--(-0.6,0.829343 )--(-0.5,0.657008 )--(-0.4,0.496973 )--(-0.3,0.350262 )--(-0.2,0.217886 )--(-0.1,0.100825 )--(0.,0)--(0.1,-0.0837514 )--(0.2,-0.149701 )--(0.3,-0.197252 )--(0.4,-0.225965 )--(0.5,-0.235566 )--(0.6,-0.225965 )--(0.7,-0.197252 )--(0.8,-0.149701 )--(0.9,-0.0837514 )--(1.,0.)--(1.1,0.100825 )--(1.2,0.217886 )--(1.3,0.350262 )--(1.4,0.496973 )--(1.5,0.657008 )--(1.6,0.829343 )--(1.7,1.01296 )--(1.8,1.20688 )--(1.9,1.41015 )--(2.,1.62186 )--(2.1,1.84117 )--(2.2,2.06729 )--(2.3,2.29948 )--(2.4,2.5371 )--(2.5,2.77952 )--(2.6,3.0262 )--(2.7,3.27665 )--(2.8,3.53043 )--(2.9,3.78713 )--(3.,4.0464 );

  \end{tikzpicture}
  \caption{A graphical description of the function $F(u)$ with only two zeros at $u=0$ and $u=\gl$.}
  \label{fig:functF}
\end{figure}
\noindent
Since $V^{\go}$ is, up to a small modification, a sum of IID variables with the same distribution as $\log \rho_1$, the function $F$ is used to compute the large deviations of $V^{\go}$, and in particular to determine the geometry of the deepest potential wells.
It is strictly convex and satisfies $F(0)=F(\gl)=0$ (see Figure \ref{fig:functF}).
We let  $u_0$ be defined by
$$F(u_0)=\min_{u\in \bbR} F(u)<0.$$

 Given a sequence of events $(A_N)_{N\ge 1}$, we say that $A_N$ holds \textit{with high probability} (which we sometimes abbreviate as w.h.p.) if $\lim_{N\to \infty}\bbP\left[A_N\right]=1$.
Given a sequence $(B_{N,k})_{N\ge 1, k\in \lint 1,N/2\rint}$, we say that $B_{N,k}$ holds 
\textit{with high probability } if 
 $$ \lim_{N\to \infty}\inf_{k \in \lint 1,N/2\rint}\bbP\left[B_{N,k}\right]=1.$$
 We are now ready to state the result.

\begin{theorem}\label{th:mixingub}
Under the assumptions \eqref{unifeliptic} and\eqref{sumpkrho}, then    with high probability we have
\begin{equation}
t_{\Mix}^{N,k,\go} \le  80 kN  \ga^{-1} \left( \frac{3 u_0+2}{\vert F(u_0)\vert}\log N \right)^4 
N^{\frac{3 u_0+2}{\vert F(u_0)\vert} \left( 2 \log \frac{1-\ga}{\ga}+4\log 4- 3\log 3 \right)}.
\end{equation}
\end{theorem}

\noindent Our second result provides a lower bound for the mixing time which depends both on $N$ and $k$.

\begin{theorem}\label{th:mixinglb}
 Under the assumptions \eqref{unifeliptic}-\eqref{sumpkrho} and assuming further that $\gl_{\bbP}<\infty$, 
there exists a positive constant $c(\ga, \bbP)$ such that 
  w.h.p. we have for every $N$ and $k\in \lint 1,N/2\rint$ 
 \begin{equation}\label{eq:mainlb}
t_{\Mix}^{N,k,\go}\geq  c \max \left\{N, N^{\frac{1}{\gl}} (\log N)^{-\frac{2}{\gl}},  k N^{\frac{1}{2\gl}} (\log N)^{-2(1+\frac{1}{\gl})}   \right \}.
\end{equation}
\end{theorem}

\subsection{Related work}\label{sec:blabla}
 Let us provide now a short review of related results present in the literature.

\subsubsection{Mixing time for the exclusion process in a homogeneous environment}
The mixing time of the exclusion process on the line segment has been extensively studied in the case where the sequence $\go$ is constant, \textit{i.e.} $\go\equiv p$. In that case, not only the right order of magnitude has been identified for the mixing time, but also the sharp asymptotic equivalent.
The case of the exclusion with no bias  that is $p=1/2$ (the simple symmetric exclusion process), it was shown  in  \cite{Aldous83} that the mixing time for the exclusion process on the segment is of order at least $N^2$ and at most $N^2(\log N)^2$. It was later established (see \cite{wilson2004mixing}
for the lower bound and \cite{lacoin2016mixing} for the upper bound) that if $k_N$ satisfies \eqref{totheinfinite}, we have 
\begin{equation}\label{determin1}
t_{\Mix}^{N,k_N}(\gep)= \frac{  (1+o(1))}{\pi^2}  N^2\log  k_N.
\end{equation}
In the case where the walk presents a bias, that is
$p\ne 1/2$, it was shown in \cite{benjamini2005mixing} that the mixing time is of order  $N$. This result was refined   
in \cite{labbe2016cutoff} by identifying the proportionality constant, showing that if $k_N$ satisfies $\lim_{N\to \infty} k_N/N=\theta$, then
\begin{equation}\label{determin2}
t_{\Mix}^{N,k_N}(\gep)=   [1+o(1)]
\frac{(\sqrt{\theta}+\sqrt{1-\theta})^2}{ |2p-1|} N.
\end{equation}
The case where $p$ is allowed to depend on $N$ was investigated in \cite{levin2016mixing, labbe2018mixing} where the order of magnitude and the sharp asymptotic of the mixing time were respectively determined. Note that in \eqref{determin1} and \eqref{determin2} the asymptotic behavior of $t_{\Mix}^{N,k_N}(\gep)$  does not display any dependence on $\gep$ at first order. This implies that $d_{N, k_N}(t)$ abruptly drops from $1$ to $0$ on the time scale $N^2 \log k_N$ and $N$ respectively. This phenomenon, called cutoff, is expected to hold for a large class of Markov chains, we refer to \cite[Chapter 18]{LPWMCMT} for an introduction. 

\medskip

Let us also mention that the mixing time for the  one-dimensional exclusion process has also been investigated for a variety of different boundary conditions. We refer to \cite{lacoinprofile} for a sharp estimate of the convergence profile to equilibrium  for the periodic boundary condition in the symmetric case and to \cite{gantertboundaries} (and references therein) for the study of a variety of boundary conditions, with or without bias. The case of higher dimension has also been considered, see e.g.\ \cite{Morris06} where the order of magnitude of the mixing time is determined up to a constant.

\subsubsection{Mixing time for the random walk in a random environment}
In \cite{gantert2012cutoff}, the case of the mixing time for a random walk in the segment with a transient random environment  (which corresponds to the case $k=1$ in the present paper) was investigated. It is shown that whenever $\gl_{\bbP}>1$ then 
\begin{equation}
 t_{\Mix}^{N,1,\go}(\gep)= [1+o(1)] \frac{N}{\bbE \left[ Q^{\go}[T^{\go}_1] \right]},
\end{equation}
where $T^{\go}_1$ is the first hitting time of $1$ for the random walk in a random environment $\go$ starting from $0$ (the result in \cite{gantert2012cutoff} is slightly more precise and the assumption is more general than \eqref{unifeliptic}). When $\gl_{\bbP} < 1$,  it is shown that the mixing time is of a much larger magnitude but that cutoff does not hold. More precisely, for $\gl_{\bbP} \le 1$ we have
\begin{equation}\label{1walk}
 \lim_{N\to \infty}  \frac{\log  t_{\Mix}^{N,1,\go}(\gep)}{ \log N}= \frac{1}{\gl_{\bbP}}.
\end{equation}
The asymptotic $N^{1/\gl_{\bbP}+o(1)}$ corresponds to the time that is required to overcome the largest potential barrier present in the system, whose height is of order $(1/\gl) \log N$.

\subsubsection{Mixing time for the exclusion in a ballistic environment}
In \cite{schmid2019mixing}, the mixing time  $t_{\Mix}^{N,k_N,\go}$ were investigated under the assumption that $\lim_{N\to \infty} k_N/N= \theta \in (0,1/2]$ and $\gl_{\bbP}>1$.
Three different cases are considered.
\begin{itemize}
 \item When $\infess \go_1>1/2$,  it is shown that the mixing $t_{\Mix}^{N,k_N,\go}$ is of order $N$, by a simple comparison with the case of homogeneous asymmetric environment.
 \item When  $\infess \go_1<1/2$, it is shown that there exists a positive $\delta$ such that the mixing time satisfies $t_{\Mix}^{N,k_N,\go}\ge N^{1+\delta}$.
 \item When $\infess \go_1=1/2$, it is shown that  
 \begin{equation}
 \liminf_{N\to \infty} t_{\Mix}^{N,k_N,\go}(\gep)/N=\infty \quad \text{ and } \quad t_{\Mix}^{N,k_N,\go}(\gep)\le C N (\log N)^3,
 \end{equation}
 together with a quantitative lower bound if $\bbP[\go_1=1/2]>0$.
\end{itemize}

 \subsubsection{Other perspectives concerning the exclusion process and  random environments}
The exclusion process with other types of random environment has also been considered in the literature.
One possibility is to consider a random environment on bonds instead of sites. A particular choice which makes the uniform measure on $\bbZ$ reversible for the random walk is the model of random conductance. In that case the mixing property of the system strongly differs from model considered here: the equilibrium measure is uniform on $\gO_{N,k}$ so that there is no trapping by potential. It is expected that for a large class of environment in that case the mixing properties are very similar to that of the homogeneous system.
The hydrodynamic limit of exclusion processes with bond–dependent  random transition rates have been studied in \cite{faggionato2008random, jara2011hydrodynamic} (see also \cite{Faggionato2020} for a recent work going slightly beyond the random conductance model). 
\medskip

\medskip

Another corpus of work has been considering the (homogeneous) exclusion process itself as a dynamical random environment, which determines the transition probability of the random walk.  
The asymptotic behavior of a random walker in this setup is studied in \cite{hilario2020random,huveneers2015random}, and the hydrodynamic limit  for the exclusion process as seen by this walker is 
studied in \cite{avena2015symmetric}. In a more general setup for the jump rates of the walker, an invariance principle about the random walk when the exclusion process starts from equilibrium is studied in \cite{jara2020symmetric}.

\subsection{Interpretation of our results and conjectures}\label{subsec:interpconject}

\subsubsection{Comments on Propositions \ref{prop:ulb} and \ref{prop:lwbdgap}}
The asymptotic for the mixing time for ASEP in homogeneous environment \eqref{determin2} shows that the lower bound of Proposition \ref{prop:ulb}
is sharp up to a constant factor.  An important observation is also that \eqref{mix2} is not true without the assumption that $k_N$ goes to infinity, even if $1/30$ is replaced by an arbitrarily small constant provided that it is not allowed to depend on $\gep$.

\medskip

However the constant in our bounds \eqref{mix1} and \eqref{mix2} are clearly not optimal.
Let us state now a natural conjecture. We believe that if 
$\lim_{N\to \infty}k_N/N=\theta\in (0,1/2]$, and $\go_x\in [\alpha,1-\alpha]$ for all $x\in \bbZ$ (with the possibility of having $\alpha=0$)
then we should have
\begin{equation}\label{laconjecture}
\liminf_{N\to \infty} \frac{1}{N}t_{\Mix}^{N,k_N}(1-\gep)\ge
\frac{(\sqrt{\theta}+\sqrt{1-\theta})^2}{1-2\alpha} . 
\end{equation}
One can obtain counter examples to \eqref{laconjecture} in the zero density case by considering the case $\go_x=1-\alpha$ in the first half of the segment $\lint 1, N\rint$ and $\go=\alpha$ in the second half of the segment, and $k_N$ diverging to infinity such that $\lim_{N\to \infty} k_N/ (\log N)=0$. In that case, one can with some minor efforts, show that the mixing time is asymptotically equivalent  $\frac{N}{2-4\alpha}$ (which is half of the lower bound in \eqref{laconjecture}). 
\medskip

Proposition \ref{prop:lwbdgap} can also be shown to be sharp within constant in the sense that there exists a constant $C_{\alpha}$, and for given $N$ and $k$ it is always possible to construct an environment $\go$ such that 
\begin{equation}
 \Gap^{\go}_{N,k} \ge e^{-C_{\alpha} N}.
\end{equation}
We conjecture that the best possible lower bound on the spectral gap  when  $\lim_{N\to \infty} k_N/ (\log N)=\theta\in (0,1/2]$ is the following
\begin{equation}
\liminf_{N\to \infty}  \frac
{\log \Gap^{\go}_{N,k_N}}{N}= -\frac{(1-\theta)}{2} \log \left(\frac{1-\alpha}{\alpha}\right).
\end{equation}
The $\liminf$ is reached asymptotically by the environment
\begin{equation}
\begin{cases}
\go_x=\alpha \quad &\text{ if } 1\le x\le \frac{(1-\theta)N}{2},\\
\go_x= 1-\alpha  & \text{ if } \frac{(1-\theta)N}{2}<  x\le N. 
\end{cases}
\end{equation}
\subsubsection{Comments on Theorems \ref{th:mixingub} and \ref{th:mixingub}}
Our paper brings a complement to the results in \cite{schmid2019mixing}, in the case when $\infess \go_1<1/2$. Firstly it provides a complementary upper bound result, which shows that the mixing time in transient environment always scales like a power of $N$, even in the non-ballistic case $\gl_{\bbP}\le 1$.

\medskip

Secondly, it provides a more quantitative lower bound. In \eqref{eq:mainlb} the mixing time is bounded by the maximum of three quantities. Each of them corresponds to a different mechanism which prevents the mixing time to be lower than a certain value.

\begin{itemize}
 \item \textit{Mass transport cannot be faster than ballistic:} Which is explored in Proposition \ref{prop:ulb} is that particle cannot move faster than ballistically (and this is independent of the choice of $\go$), so that the time required to transport the mass of particles to equilibrium has to be at least of order $N$. This idea is already present in \cite{benjamini2005mixing}.
 \item \textit{Individual particles may be blocked by traps in the potential profile:}
 As soon as $\infess \go_1<1/2$, the potential profile $V$ is non-monotone and will create energy barriers.  It is  known since \cite{kesten1975limit}  that these energy barriers can slow down particles to subballistic speed in $\gl_{\bbP}\le 1$ by creating traps that will require a long time to be crossed. This is the mechanism that was used to identify the mixing time in case of a single particle in \cite{gantert2012cutoff} (recall \eqref{1walk}), and it corresponds to the time needed to cross the largest trap in the potential. This yields the second term in  \eqref{eq:mainlb}.
 
 \item \textit{Potential barrier may also create bottleneck for the flow of particles:}
 The third mechanism which was partially identified in \cite{schmid2019mixing} is that potential barrier may also limit the flow of particles throughout the system.
 The limitation on the flow does not correspond to the inverse of the time that a particle needs to cross the trap, but rather to the square root of this inverse.
 The reason for this is that when particles are flowing through the system, the particle are ``filling'' half of the potential well, so that the remaining potential barrier to be crossed is halved.
 This reasoning yields the third term in \eqref{eq:mainlb}.
\end{itemize}

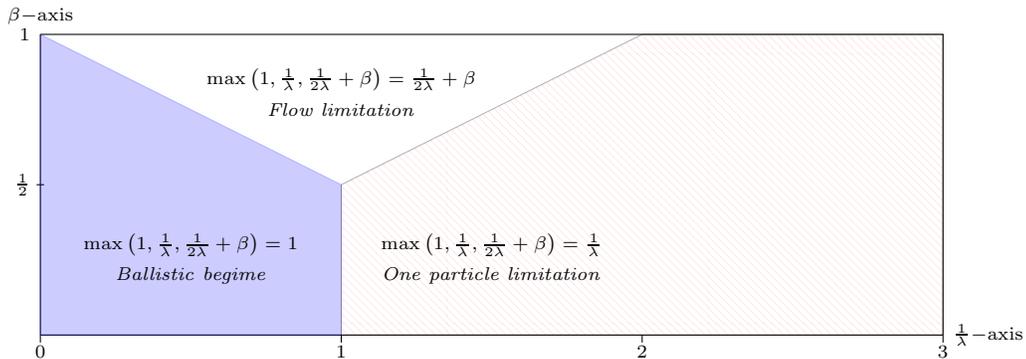
\begin{figure}[h]
  \centering  
\begin{tikzpicture}[scale=1,font=\tiny] 
\draw (0,0)--(12,0);
       \draw (0,0) -- (0,4);       
      \foreach \x in {0,4,8,12} {  \draw (\x,0)--(\x,-0.1);}
     \draw (12,0) -- (12,4)--(0,4);
     \node[below] at (0, 0){$0$};
       \node[below] at (4,0){$1$};
      \node[below] at (8,0){$2$};
      \node[below] at (12,0){$3$};
      \node[left] at (0,4){$1$};
      \node[above] at (0,4){$\gb-$axis};
      \node[right] at (12,0){$\frac{1}{\gl}-$axis};
      \node[left] at (0,2){$\frac{1}{2}$};
      \draw(-0.05,2)--(0.05,2);

\draw[color=blue,fill,opacity=.2] (0,4)--(4,2)--(4,0)--(0,0)--(0,4);
    
 \draw[pattern=north west lines, pattern color=red, opacity=.3] (4,0)--(4,2)--(8,4)--(12,4)--(12,0)--(4,0);  

  \node at (2,1.2){$\max\left(1, \frac{1}{\gl}, \frac{1}{2\gl}+\beta \right)=1$};
\node at (2,0.8) {\textit{Ballistic begime}}; 
  \node at (6,1.2){$\max\left(1, \frac{1}{\gl}, \frac{1}{2\gl}+\beta \right)=\frac{1}{\gl}$};
    \node at (6,0.8){\textit{One particle limitation}};
  \node at (4,3.4){$\max\left(1, \frac{1}{\gl}, \frac{1}{2\gl}+\beta \right)=\frac{1}{2\gl}+\gb$};
    \node at (4,3.0) {\textit{Flow limitation}};
  \end{tikzpicture}
  \caption{The phase diagram for the exponent of the mixing time (the lower bound is proved rigorously and the upper bound is only conjectured). The transition between the blue and red (hatched) regions of the diagram corresponds to the transition of the RWRE from the ballistic phase to the transient-with-zero-speed phase. A third phase represented by the white region appears when one considers a large number of particles, in this phase the main limitation to mixing is the flow of particle through the deepest trap.
  } \label{fig:dynamicdiagram}
\end{figure}

We believe that the three mechanism described above are the only possible limiting factor to mixing, and thus that the lower bound give in Theorem \ref{th:mixinglb} is sharp as far as the exponent is concerned.
Let us formulate this as a conjecture.
Let us assume that $k_N$ satisfies $$\lim_{N\to \infty} \frac{\log k_N}{\log N}=\beta,$$
and then we should have the following convergence w.h.p.
\begin{equation}
 \lim_{N\to \infty} \frac{t_{\Mix}^{N,k_N}}{\log N}= \max\left(1, \frac{1}{\gl}, \frac{1}{2\gl}+\beta \right).
\end{equation}
We refer to Figure \ref{fig:dynamicdiagram} for the phase diagram concerning the conjectured exponent of the mixing time.

In particular this means that when $\beta \le  1/(2\gl)$ then the mixing time of the exclusion process on the segment coincides (as far as the exponent is concerned) with that of the random walk in the segment.

\subsection*{Organization}  

Section \ref{sec:techpre} is devoted to some technical preliminaries including the particle description, equilibrium estimates, partial order, a graphical construction and a composed censoring inequality.

Section \ref{sec:univbounds} is devoted to  universal lower and upper bounds on the mixing time for all random environments, that is, the proof of Propositions \ref{prop:ulb} and Proposition \ref{prop:lwbdgap}.

Section \ref{sec:lowbdmix} is devoted to lower bounds on the mixing time, that is Theorem \ref{th:mixinglb}. There are three bounds to prove, one of them is a consequence of Proposition \ref{prop:ulb}, the other two are presented as  two distinct results (Proposition \ref{prop:oneparticle} and Proposition \ref{prop:flow}) and proved in separate subsections. The first bound rely on controlling the displacement of the leftmost particle while the other is based on a control of the particle flow.

Section \ref{sec:upbdmix} is concerned with the upper bound on the mixing time (Theorem \ref{th:mixingub}). The proof is based on application of the censoring inequality and of our upper bound from Proposition \ref{prop:lwbdgap}: blocking the transition along carefully chosen edges (in a way that varies throught time) we guide all particlesto the right of the segment (where they are typically located at equilibrium) in polynomial time.

\subsection*{Notation} 
We use $c(\ga, \bbP)$ and $C(\ga,\bbP)$ to stress that the constants $c$ and $C$ depend on $\ga$ and the law of the random environment $\go$. Moreover, we use $``\colonequals"$ (or $``\equalscolon"$) to define a new  quantity on the left-hand (right-hand, resp.) side, and $\lint a, b\rint \colonequals [a,b] \cap \bbZ$. Furthermore, we let
\begin{equation}
\gO_{[a,b],k} \colonequals \left\{ \xi \in \{0,1 \}^{\lint a, b \rint}: \sum_{x \in \lint a, b \rint}\xi(x)=k  \right\}
\end{equation}
denote the state space of $k$ particles performing exclusion process restricted in the interval $\lint a, b\rint$ and environment $\go$, and let $\pi_{[a,b],k}^{\go}$ denote the corresponding  equilibrium probability measure.

\subsection*{Acknowledgment}
The authors thank Milton Jara, Roberto Imbuzeiro Oliveira and  Augusto Teixeira for enlightening discussions. This work was realized in part during H.L.\ extended stay in Aix-Marseille University funded by the European Union’s Horizon 2020 research and innovation programme under the Marie Skłodowska-Curie grant agreement No 837793.

\section{Technical preliminaries}\label{sec:techpre}

 \subsection{Partial order on $\gO_{N,k}$}

Given $\xi \in \gO_{N,k}$ we define  $\bar \xi \ : \  \lint 1, k\rint \to \lint 1, N\rint$ as an increasing function which provides the  positions of the particles of $\xi$ from left to right:

\begin{equation}\label{partiposition}
 \{\bar \xi(i)=x\}  \quad  \Longleftrightarrow  \quad  \left\{\xi(x)=1 \text{ and } \sum_{y=1}^x \xi(x)=i \right\}.
\end{equation}
We introduce a natural partial order relation ``$\leq$'' on $\Omega_{N,k} \times \Omega_{N,k}$  as follows
\begin{equation}\label{partorder}
\left( \xi \leq \eta \right)  \mbox{ } \Leftrightarrow \left(\forall i \in \lint 1, k\rint,  \quad \bar \xi(i)\leq \bar \eta (i) \right). 
\end{equation}
Informally $\xi \leq \eta$  means that the particles in the  configuration $\eta$ are located ``more to the right'' than those of $\xi$. 
Let $\xi_{\max}$ and $\xi_{\min}$ denote the maximal and minimal configurations of $(\Omega_{N,k}, ``\leq")$ respectively, given by
\begin{equation}\label{minmaxconf}
\begin{aligned}
 \xi_{\max} \colonequals \ind_{\{N-k+1 \le x\le N\}} \quad \mbox{ and } \quad
  \xi_{\min} \colonequals \ind_{\{1\leq x \leq k\}}.
\end{aligned}
\end{equation}
This order plays a special role for our dynamic $(\sigma^{\xi}_t)_{t \ge 0}$, and the next two subsections provide tools to exploit this link.

\subsection{Canonical coupling via graphical construction}\label{subsec:cancoupling} 
Let us present a construction of a grand coupling for the exclusion process on the segment $\lint 1,N\rint$ which has the property of conserving the order defined above.

\medskip

To each site $x \in \lint 1, N \rint$ we associate an independent rate $1$ Poisson clock process $(T^{(x)}_i)_{i\ge 1}$ (the increments of the sequence $(T^{(x)}_i)_{i \ge 1}$ are independent exponential variables of parameter $1$) and an independent sequence of IID variables $(U^{(x)}_i)_{i\ge 1}$ with uniform distribution on $[0,1]$.
These variables are independent of the environment $\omega=(\go_x)_{x \in \bbZ}$, and 
the trajectory $(\sigma^\xi_t)_{t\ge 0}$ for each $\xi$ is a deterministic function of $(T^{(x)}_i, U^{(x)}_i)_{i\ge 1, x\in \lint 1,N\rint }$.
In the remainder of the paper, $\bP$ denote the joint law of $(T^{(x)}_i, U^{(x)}_i)_{i\ge 1, x\in \lint 1,N\rint }$, and $\bE$ denotes the corresponding expectation.
Let us also introduce a natural filtration $(\cF_t)_{t\ge 0}$ in this probability space 
setting 
\begin{equation}\label{defizero}
i_0(x,t):=\max\{ i\ge 1 \ : \ T^{(x)}_i\le t\}
\end{equation}
with the convention that $\max \emptyset=0$ and set
\begin{equation}\label{filtrations}
 \cF_t:= \sigma\left( T^{(x)}_i, U^{(x)}_i, x\in \bbZ, i\le i_0(x,t)  \right).
\end{equation}
Now, given $1\le k\le N-1$ and an initial configuration $\xi\in \gO_{N,k}$, we construct the trajectory 
$(\sigma^{\xi}_{t})_{t \ge 0}$ 
as follows:
\begin{itemize}
\item[(1)] $(\sigma^{\xi}_t)_{t \ge 0}$ is càdlàg and may change its value only  at times $T^{(x)}_i$, $x\in \lint 1,N\rint$ and $i\ge 1$.
\item[(2)] We construct the trajectory starting with $\sigma^{\xi}_0=\xi$ and modifying it sequentially  at the update times 
 $(T^{(x)}_i)_{i\ge 1, x\in \lint 1,N\rint}$. For instance if $t=T^{(x)}_i$ we obtain $\sigma^
 {\xi}_{t-}$ from $\sigma^
 {\xi}_{t}$ as follows:
\begin{itemize}
\item[(A)] If $U^{(x)}_i\le \go_x$, $ x \le N-1$,  $\sigma^{\xi}_{t_-}(x)=1$ and $\sigma^{\xi}_{t_-}(x+1)=0$,
then $\sigma^{\xi}_{t}(x+1)=1$ and $\sigma^{\xi}_{t}(x)=0$ (and $\sigma^{\xi}_{t}(y)=\sigma^{\xi}_{t_-}(y)$ for $y\notin \{x,x+1\}$).
\item[(B)] If $U^{(x)}_i> \go_x$, $ x \ge 2$, $\sigma^{\xi}_{t_-}(x)=1$ and $\sigma^{\xi}_{t_-}(x-1)=0$,
then $\sigma^{\xi}_{t}(x-1)=1$ and $\sigma^{\xi}_{t}(x)=0$  (and $\sigma^{\xi}_{t}(y)=\sigma^{\xi}_{t_-}(y)$ for $y\notin \{x-1,x\}$).
\item[(C)] In all other cases $\sigma^{\xi}_{t}=\sigma^{\xi}_{t_-}$.
\end{itemize}
\end{itemize}

It is immediate by inspection to check that the above construction corresponds indeed to the Markov chain with generator $\cL^{\go}_{N,k}$.
Note also that  our process is adapted and Markov with respect to the filtration $(\cF_t)_{t \ge 0}$.
 In the same manner, the reader can check that it preserves the order in the following sense.

\begin{proposition}\label{monocouple}
 For the coupling constructed above, we have for all $\xi, \xi' \in \gO_{N,k}$
 
 \begin{equation}
\xi\leq \xi'  \Rightarrow  \bP\left[ \forall t\geq 0, \quad \sigma^{\xi}_t\leq \sigma_t^{\xi'}   \right]=1.
\end{equation}
\end{proposition}

\subsection{Composed censoring inequality}\label{compcensor}

We are going to use a variant of the censoring inequality introduced by Peres and Winckler 
\cite{peres2013can}. Let $E_N=\left\{ \left\{n,n+1 \right\}: n\in \lint 1, N-1\rint \right\}$ be the set of edges in $\lint 1, N \rint$, and a censoring scheme $\mathcal{C}: [0, \infty) \to \mathcal{P}(E_N)$ is a deterministic c\`{a}dl\`{a}g function where $\mathcal{P}(E_N)$ is the set of all subsets of $E_N$.

\medskip

The censored chain $(\sigma^{\xi,\cC}_t)_{t \ge 0}$ is a time inhomegenous Markov chain, with a generator
obtained by cancelling the transition using edges in $\cC(t)$

\begin{equation}\label{cgenerator}
\mathcal{L}^{\cC,t}_{N,k}(f)(\xi) \colonequals \sum_{x=1}^{N-1} r^{\go}_{N,k}(\xi,\xi^{x,x+1}) \ind_{\left\{\{x,x+1\}\notin \cC(t)\right\}}
\big[f(\xi^{x,x+1})-f(\xi)\big],
\end{equation}
where $r^{\go}_{N,k}(\xi,\xi^{x,x+1}) $ is defined in \eqref{jumprateSEP}.
We let  $P^{\cC}_t$ be the associated semigroup (the solution of $\partial_t P_t= P_t\cL_{N,k}^{\cC,t}$ with initial condition given by the identity).
We will use the following corollary of the censoring inequality \cite[Theorem 1]{peres2013can} (recall \eqref{minmaxconf}).
\begin{proposition}\label{PWcoro}
 For any $\xi\in \gO_{N,k}$ and any censoring scheme $\cC$, we have 
 \begin{equation}
  P_t(\xi,\xi_{\max})\ge P^{\cC}_t(\xi_{\min},\xi_{\max})
 \end{equation}

\end{proposition}

\begin{proof}[Sketch of proof]
 Proposition \ref{monocouple} implies that  $P_t(\xi,\xi_{\max})\ge P_t(\xi_{\min} ,\xi_{\max})$.
 To compare $P_t(\xi_{\min} ,\xi_{\max})$ with  $P^{\cC}_t(\xi_{\min} ,\xi_{\max})$, we  rely on the censoring inequality \cite[Theorem 1]{peres2013can} (to see that the exclusion process fits the setup in \cite{peres2013can}, one uses the height function representation see e.g.\ \cite[Section A.2]{lacoin2016mixing}) which implies that $P_t(\xi_{\min} ,\cdot)$ stochastically dominates $P_t^{\cC}(\xi_{\min} ,\cdot)$.
 
\end{proof}

\medskip

We consider a modified censored dynamics, where on top of censoring,  at fixed time, we replace the current configuration by one which is lower for the order $\ge$ by moving some particles to the left. For the application we have in mind, we can consider that these replacements are performed deterministically (although the result would hold also for random replacements).

\medskip

Let $(s_i)_{i=1}^I$ be an increasing time sequence tending to infinity and let  $(Q_i)_{i=1}^I$ be a sequence of stochastic matrices on $\gO_{N,k}$ such that for all $\xi$ in $\gO_{N,k}$ there exists 
$\xi'$ (depending on $\xi$ and $i$) such that 
\begin{equation}\label{upp}
\begin{cases}
 \xi' \le \xi,\\
 Q_i(\xi, \xi')=1,\\
 Q_i(\xi,\xi'')=0, \text{ when } \xi''\ne \xi'.
\end{cases}
\end{equation}
We consider $\tilde P_t$ the semigroup defined by
\begin{equation}\label{semigroup:displace}
 \begin{cases}
\tilde P_0= \mathrm{Id},\\
\partial_t \tilde P_t= \tilde  P_t \cL^{\cC,t} \text{ if } t\notin \{s_i\}_{i=1}^I,\\
\tilde P_{s_i}=\tilde P_{(s_i)_-} Q_i.
 \end{cases}
\end{equation}

\begin{proposition}\label{prop:censor}
For any choice of $(s_i)_{i=1}^I$, $(Q_i)_{i=1}^I$ and $\cC$, we have  for all $t \ge 0$
 \begin{equation}
   P^{\cC}_t(\xi_{\min},\xi_{\max}) \ge \tilde P_t(\xi_{\min},\xi_{\max}).  
\end{equation} 
\end{proposition}

\begin{proof}
 We construct both $(\tilde \sigma^{\min}_t)_{t \ge 0}$ with transition probability $\tilde P_t$ with initial condition $\xi_{\min}$ and $(\sigma^{\min, \cC}_t)_{t \ge 0}$  the censored dynamics with the same initial condition on the same probability space, using the variables $(T^{(x)}_i, U^{(x)}_i)_{i\ge 1, x\in \lint 1,N\rint }$.
 
 \medskip
 
For $(\sigma^{\min, \cC}_t)_{t \ge 0}$ we use the same procedure  as for $(\sigma^{\xi}_t)_{t\ge 0}$ (for $\xi=\xi_{\min}$) with the following added requirement for the transitions: $\{x,x+1\}\notin \cC(t)$ in the case $(A)$ and $\{x,x-1\}\notin \cC(t)$ in the case $(B)$. \\
 For  $(\tilde \sigma^{\min}_{t})_{t\ge 0}$ we use the same procedure as  for $(\sigma^{\min, \cC}_t)_{t \ge 0}$ but with the addition of new deterministic jumps in the trajectories at times $(s_i)_{i\in I}$. 
More precisely if  $t=s_i$,  $\tilde \sigma^{\min}_{t}$ is determined from  $\tilde \sigma^{\min}_{t_-}$ as the unique element of $\gO_{N,k}$ such that
\begin{equation}
 Q_i (\tilde \sigma^{\min}_{t_-},\tilde \sigma^{\min}_{t})=1.
\end{equation}
We have by definition $\tilde \sigma^{\min}_0=\sigma^{\min, \cC}_0$, and it can be checked by inspection that all the transitions are order preserving (this is a property of the graphical construction when $t\notin \{s_i\}_{i=1}^I$ and a consequence of \eqref{upp} for  the special values $t\in \{s_i\}_{i=1}^I$).

\end{proof}

\subsection{Equilibrium estimates}\label{subsec:equilesti}

Recalling \eqref{loglaplace}
let us define 
\begin{equation}
 \kappa:= F'(\gl)=\bbE\left[\rho_1^{\gl} \log (\rho_1) \right]>0,  
\end{equation}
and set
\begin{equation}\label{deltamax}
 \Delta V^{\go,N}_{\max}= \max_{ 1\le x \le y \le N} \left( V(y)-V(x) \right).
\end{equation}
The literature on the subject of random walks in a random environment contains very sharp information concerning  $\Delta V^{\go,N}_{\max}$, and the length of the corresponding trap (see  \cite{gantert2012cutoff}). In particular it is known under quite general assumptions that
$|\Delta V^{\go,N}_{\max}- \frac{1}{\gl}\log N|$ displays random fluctuations of order $1$
and that the corresponding traps are of a length $\frac{1}{\gl \kappa}\log N$ at first order.

\medskip

For the sake of completeness we include a short proof of the following non-optimal result which is sufficient to our purpose.
 Set
\begin{equation}\label{def:qN}
q_N \colonequals \frac{3u_{0}+2 }{|F(u_{0})|} \log N,
\end{equation}
where $u_0$ is the point at which $F$ attains its minimum.

\begin{proposition}\label{lema:createtrap}
We have
\begin{equation}\label{deltavalue}
 \lim_{N\to \infty}\bbP\left[  -\left(\frac{1+\gep} \gl \right)\log \log N  \le \Delta V^{\go,N}_{\max}- \frac{1}{\gl}\log N\le  \frac{\gep}{\gl} \log \log N  \right]=1.
\end{equation}
Furthermore we have
\begin{equation}\label{mvalue}
 \lim_{N\to \infty}\bbP\left[ \maxtwo{ 1\le x \le y \le N}{y-x\ge q_N} \left( V(y)-V(x) \right) \ge -3\log N\right]=0.
\end{equation}
In particular, with high probability we have 
$$\forall x,y\in \lint 1,N\rint, \mbox{  } \left\{ V(y)-V(x)= \Delta V^{\go,N}_{\max} \right\} \quad  \Rightarrow \left\{(y-x)\le q_N\right\}.$$

\end{proposition}

\begin{proof}

At the cost of an additive constant on our bounds (which we omit in the proof for readability), using our uniform ellipticity assumption we can replace 
$V(y)-V(x)$ in the definition of \eqref{deltamax}
by a sum of IID random variables,  setting $\bar V(1)=0$ and
\begin{equation}\label{booz}
\sum_{z=x+1}^{y} \log \rho_z:= \bar V(y)- \bar V(x).
\end{equation}
By definition  of $\gl$, $M_n=\left(\prod_{x= 1}^n (\rho_x)^{\gl}\right)_{n\ge 1}$ is a martingale  for the filtration $\cG_n:= \sigma(\go_x, x\in \lint 1,n\rint)$.
Using the optional stopping theorem at  $T_A:=\inf\{ n, M_n\ge A\}$ and using that 
\begin{equation}
\begin{cases}
 A\le  M_{T_A}\le A\left(\frac{1-\alpha}{\alpha}\right)^{\gl},\\
\lim_{n\to \infty} M_n=0,
  \end{cases}
\end{equation}
we have for any $A$
\begin{equation}\label{foranyA}
\frac{1}{A} \left(\frac{\alpha}{1-\alpha}\right)^{\gl}\le  \bbP\left[ \max_{n\ge 1}\prod_{x= 1}^n (\rho_x)^{\gl} \ge A \right]\le \frac{1}{A}.
\end{equation}
The bound above can be used to obtain the upper bound on $\Delta V^{\go,N}_{\max}$ via a union bound using translation invariance
\begin{multline}
 \bbP\left[  \max_{1\le x\le  y \le N}\bar V(y)-\bar V(x) \ge \frac{1}{\gl}\log N + \frac{\gep}{\gl} \log \log N\right]\\
 \le \sum_{x=1}^{N} 
 \bbP\left[ \max_{y\ge x } \bar V(y)-\bar V(x) \ge  \frac{1}{\gl}\log N + \frac{\gep}{\gl} \log \log N  \right]\\
 \le N 
 \bbP\left[ \max_{n\ge 1}\prod_{x= 1}^n (\rho_x)^{\gl} \ge N (\log N)^{ \gep} \right]
 \le (\log N)^{- \gep}.
\end{multline}
Before proving the corresponding lower bound, let us move to the proof of \eqref{mvalue}. Again using translation invariance  and union bound, it is sufficient to show that 
\begin{equation}
 \lim_{N\to \infty} N
    \bbP\left[ \max_{n \ge q_N }\sum_{x=1}^n \log \rho_x \ge -3\log N \right]=0.
\end{equation}
We use Doob's maximal inequality for the martingale  $e^{-nF(u_0)}\prod_{x=1}^n (\rho_x)^{u_0}$.
Since $F(u_0)< 0$, we have
 \begin{multline}\label{zbad}
   \bbP\left[ \max_{n \ge q_N }\prod_{x=1}^n (\rho_x)^{u_0} \ge  N^{-3u_0} \right]
   \le    \bbP\left[ \max_{n \ge 1 } e^{-nF(u_0)}\prod_{x=1}^n (\rho_x)^{u_0} \ge  N^{-3 u_0} e^{-q_N F(u_0)} \right]
  \\ \le N^{3u_0}e^{q_N F(u_0)}\le N^{-2}.
 \end{multline}
 This is sufficient to conclude the proof of \eqref{mvalue}.
 Note that as a consequence by \eqref{foranyA} and \eqref{zbad}, we have for $N$ sufficiently large
 \begin{equation}
     \bbP\left[ \max_{1\le n \le q_N }\prod_{x=1}^n (\rho_x)^{\gl} \ge  N (\log N)^{-(1+\gep)} \right]
      \ge \frac{1}{2}\left(\frac{\alpha}{1-\alpha}\right)^{\gl} N^{-1}(\log N)^{1+\gep}.
 \end{equation}
As a consequence of independence we have
\begin{multline}
 \bbP \left[ \forall (i,j) \in \lint 1, \lfloor N/q_N \rfloor-1 \rint\times \lint 1,q_N\rint, \   \bar V(iq_N+j)-\bar V(iq_N)\le  \frac{\log N- (1+\gep)\log \log N}{\gl}\right]\\ \le \left(1- \frac{1}{2}\left(\frac{\alpha}{1-\alpha}\right)^{\gl}N^{-1}(\log N)^{(1+\gep)}\right)^{\lfloor N/q_N \rfloor-1} \le e^{-c (\log N)^{\gep}}
\end{multline}
This yields the lower bound in \eqref{deltavalue}.

\end{proof}

\begin{proof}[Proof of Lemma \ref{lema:equilleftright}]
For $\xi \in \gO_{N,k}$,   we define the positions of its leftmost particle and rightmost empty site to be respectively
\begin{equation}
\begin{gathered}
L_{N,k}(\xi) \colonequals \inf \left\{ x \in \lint 1, N \rint: \ \xi(x)=1   \right\},\\
R_{N,k}(\xi) \colonequals \sup \left\{ x \in \lint 1, N \rint: \ \xi(x)=0   \right\}.
\end{gathered}
\end{equation}
Then $$ \pi_{N,k}^{\go} \left(\cA_r^{\complement} \right) \le \pi_{N,k}^{\go}  \left(  L_{N,k}(\xi) \le N-k-r \right)+ \pi_{N,k}^{\go} \left( R_{N,k}(\xi) \ge N-k+r \right).$$
 Let us bound the second term, the first one can be treated in a symmetric manner. Moreover, we have
\begin{equation}
\pi_{N,k}^{\go} \left( R_{N,k}(\xi) \ge N-k+r \right)= \sumtwo{x \in  \lint 1, N-k\rint}{y \in \lint N-k+r,N\rint} \pi_{N,k}^{\go} \left( L_{N,k}=x, R_{N,k}=y \right)
\end{equation}
Furthermore, we recall that $\xi^{x,y}$,  defined in \eqref{swap2sites}, denotes the configuration obtained by swapping the values at sites $x,y$ of the configuration $\xi$, and observe that the map $\xi\to \xi^{x,y}$ is injective  from $\{ \xi \in \gO_{N,k}: L_{N,k}(\xi)=x, R_{N,k}(\xi)=y \}$ to $\gO_{N,k}$ defined by $\xi \mapsto \xi^{x,y}$. Then we have  
\begin{multline}\label{exchange01}
\pi_{N,k}^{\go} \left( L_{N,k}=x, R_{N,k}=y \right)
=\sum_{  \{\xi: \ L_{N,k}(\xi)=x, R_{N,k}(\xi)=y \}} \pi_{N,k}^{\go}(\xi^{x,y}) e^{V^{\go}(y)-V^{\go}(x)}\\ \le e^{V^{\go}(y)-V^{\go}(x)}
\le C e^{\bar V^{\go}(y)-\bar V^{\go}(x)}.
\end{multline}
Now from the law of large number applied to sum of IID variables, we have
\begin{equation}
 \lim_{r\to \infty} \inftwo{N\ge 1}{ k\in \lint 1, N/2\rint}\bbP\left[ \forall(x, y) \in  \lint 1, N-k\rint\times \lint N-k+r,N\rint, \  \bar V^{\go}(y)-\bar V^{\go}(x)\le \frac{ (y-x)  \bbE[\log \rho_1]}{2}\right]=1.
\end{equation}
Moreover, since $$\sumtwo{x \in  \lint 1, N-k\rint}{y \in \lint N-k+r,N\rint} e^{\frac{ \bbE[\log \rho_1] (y-x)}{2}}\le \frac{ e^{ \bbE[\log \rho_1] r/2}}{(1-e^{ \bbE[\log \rho_1]/2})^2}$$
 we have 
\begin{equation}\label{prob:rempty}
  \lim_{r\to \infty} \inftwo{N\ge 1}{ k\in \lint 1, N/2\rint}\bbP\left[ \pi_{N,k}^{\go} \left( R_{N,k}(\xi) \ge N-k+r \right) \le  \left(1-e^{\frac { \bbE[\log \rho_1]} 2}\right)^{-2}e^{\frac{{ \bbE[\log \rho_1]} r}2} \right]=1,
\end{equation}
which concludes the proof.

\end{proof}

\section{Bounds for the mixing time with arbitrary environments}\label{sec:univbounds}

\subsection{Proof of Proposition  \ref{prop:ulb}}

We look at the variable 
$$ m(\xi):= \sum_{x=1}^N x \xi(x).$$
Note that $m(\xi)\in \left[ \frac{k(k+1)}{2}, \frac{k(2N-k+1)}{2} \right]$.
We assume that
$$\pi^{\go}_{N,k}\left(m(\xi)\ge \frac{k(N+1)}{2}\right)\ge 1/2$$
(the other case can be treated symmetrically).
Now, since at all time each particle jumps to right with a rate which is at most one, starting from $\xi_{\min}$  (we write $\sigma^{\min}_t$ for $\sigma^{\xi_{\min}}_t$ to lighten the notation) we have 
\begin{equation}
 \bE\left[ m(\gs^{\min}_t)\right]\le \frac{k(k+1)}{2}+ kt.
\end{equation}
As a consequence of Markov's inequality, we have
\begin{equation}
 \bP\left[m(\gs^{\min}_t) \ge \frac{k(N+1)}{2} \right]=   \bP\left[m(\gs^{\min}_t)-  \frac{k(k+1)}{2}\ge \frac{k(N-k)}{2} \right]\le \frac{2t}{(N-k)},
\end{equation}
which is smaller than $1/4$ if $t\le N/16$.

\medskip

When the number of particles goes to infinity, we use the same kind of reasoning but adding concentration estimates for $m(\xi)$, under the equilibrium measure $\pi^{\go}_{N,k}$ (which is denoted simply by $\pi$ in this proof for readability).
Let us prove that 
\begin{equation}\label{varvar}
\Var_{\pi}[m(\xi)]\le N^2 k.
\end{equation}
To this end we introduce the filtration $(\cG_i)^N_{i=1}$ defined by $\cG_i:= \sigma(\xi(x), x\in\lint 1,i\rint)$, and consider the martingale
\begin{equation}
 M_i:= E_{\pi}\left[m(\xi) \ | \ \cG_i\right]
\end{equation}where $E_{\pi}[ \cdot \  | \ \cG_i]$ denotes the conditional expectation under $\pi$ .
We have by construction
\begin{equation}
 \Var_{\pi}[m(\xi)]= \sum_{i=1}^N \Var(M_i-M_{i-1})
\end{equation}
Now, we are going to show that 
\begin{equation}\label{swoob}
 \Var(M_i-M_{i-1})\le  \pi(\xi_i=1) (N-i)^2
\end{equation}
which implies \eqref{varvar}.
To prove \eqref{swoob} we are going to show that for any $\chi\in \{0,1\}^{i-1}$ with at most $k-1$ ones and at most $N-k-1$ zeros, the quantity
\begin{equation}\label{sbob}
 \Delta_i(\chi)= E_{\pi}\left[m(\xi) \ |\ \xi_{ \lint 1,i-1\rint}= \chi, \xi(i)=0 \right]- 
 E_{\pi}\left[m(\xi) \ | \ \xi_{\lint 1,i-1\rint}=\chi, \xi(i)=1 \right]
 \end{equation}
 satisfies
 \begin{equation}
   0 \le  \Delta_i(\chi) \le N-i.
 \end{equation}
Note that we have
\begin{equation}\label{equilmeasrestrict}
 E_{\pi}\left[m(\xi) \ |\ \xi_{ \lint 1,i-1\rint}= \chi \right]= \sum_{x=1}^{i-1} x \chi(x) + \pi^{\go}_{\lint i,N \rint,k- \sum_{x=1}^{i-1} \chi (x)} \left( \sum_{x=i}^N x \xi(x)\right),
\end{equation}
where if $I$ is a segment on $\bbZ$ and $k'\le |I|$, $ \pi^{\go}_{I,k'}$ denotes the equilibrium measure for exclusion process on $I$ with $k'$ particles and environment $\go$.
For this reason it is sufficient to prove \eqref{sbob} for $i=1$,  and arbitrary $k$ (not necessarily assuming $k\le N/2$).
Hence we need to prove that for  $N\ge 1$ and $k\in \lint 1, N-1\rint$ we have
\begin{equation}\label{trop1}
0 \le E_{\pi}\left[m(\xi) \ |\ \xi(1)=0\right]- E_{\pi}\left[m(\xi) \ |\ \xi(1)=1\right] \le N-1.
\end{equation}
 To prove this we observe that there exists a probability $\Pi$ on $\gO^2_{N,k}$  with marginals $\pi(\cdot \ |\ \xi(1)=0)$ and $\pi(\cdot \ |\ \xi(1)=1)$ such that
\begin{equation}\label{trop}
 \Pi\left( \sum_{x=1}^{N} \ind_{\{\xi^1(x)\ne \xi^2(x)\}}=2\right)=1
\end{equation}
(meaning that $\xi^1(x)=\xi^2(x)$ except at two sites, $1$ and another random site).
With this coupling we have 
\begin{equation}
 E_{\pi}\left[m(\xi) \ |\ \xi(1)=0\right]- E_{\pi}\left[m(\xi) \ |\ \xi(1)=1\right]
 = \Pi\left[ \sum_{x=1}^N x(\xi^1(x)-\xi^2(x))\right],
\end{equation}
which yields \eqref{trop1}.
The coupling $\Pi$ can be achieved using the graphical construction: we define $(\xi^1_t)$ and $(\xi^2_t)$ starting  with initial configuration $\ind_{\lint 2,k+1\rint}$ and $\ind_{\lint 1,k\rint}$ respectively and evolving using the graphical construction with the edge $\{1,2\}$ censored (recall Section \ref{compcensor}). The dynamic conserves the number of discrepancy and 
 $\pi(\cdot |\ \xi(1)=0)$ and $\pi(\cdot |\ \xi(1)=1)$ are the respective equilibrium distribution of the marginals, so that any limit point of $\bP\left[ (\xi^1_t,\xi^2_t)\in \cdot  \right]$ (existence is ensured by compactness) provides a coupling satisfying \eqref{trop}.
 
 \medskip

\noindent Now to see that \eqref{sbob} implies \eqref{swoob}, we simply observe that, 
conditionned to the state of the first $i=1$ vertices of the segment, $(M_i-M_{i-1})$ can only assume two values which differ by an amount    $\Delta_i(\chi)$ (cf. \eqref{sbob}). The corresponding conditioned variance is equal to $\Delta_i(\chi)^2$ times that of the corresponding Bernoulli variable that is
\begin{multline}
 E_{\pi}[ (M_i-M_{i-1})^2 \ | \ \xi_{ \lint 1,i-1\rint}= \chi ]= \pi ( \xi(i)=1 \ |\ \xi_{ \lint 1,i-1\rint}= \chi ) \pi(\xi(i)=0 \ | \ \xi_{ \lint 1,i-1\rint}= \chi)
 \Delta_i(\chi)^2 \\ 
 \le  \pi ( \xi(i)=1 \ |\ \xi_{ \lint 1,i-1\rint}= \chi ) (N-i)^2.
\end{multline}
We then consider the average the inequality with respect to $\xi_{ \lint 1,i-1\rint}$ to conclude.
Now using \eqref{varvar} we can assume that for any $\gep$ there exists $N_0(\gep)$ such that  for $N\ge N_0(\gep)$ we have 
\begin{equation}
\min\left[ \pi^{\go}_{N,k}(m(\xi)\le N k /3), \pi^{\go}_{N,k}(m(\xi)\ge 2Nk/3) \right]\le \gep/2.
\end{equation}
Let us assume  that the first of these two terms is smaller (the other case is treated symmetrically). To conclude, we must show that for $t=\frac{N}{30}$  we have
\begin{equation}
 \bP\left(m(\gs^{\min}_t)> N k /3\right)\le  \gep/2. 
\end{equation}
To check this we observe that 
\begin{equation}
m(\gs^{\min}_t)\le \frac{k(k+1)}{2}+\cN_t
\end{equation}
where $\cN_t$ is the total number of particle jumps to the right. Since each particle jumps at most with rate one, we have for $N$ sufficiently large 
\begin{equation}
 \bP\left[ \cN_t \ge 2kt \right] \le \gep/2,
\end{equation}
which allows to conclude. 

\qed

\subsection{Proof of Proposition  \ref{prop:lwbdgap}}

 For the proof of Proposition \ref{prop:lwbdgap}, we apply the so-called flow method (see \cite[Chapter 13.4]{LPWMCMT}).
  A path  $\gG$ is a sequence of configurations
 $(\xi_0,\dots, \xi_{|\gG|})$ which is such that  $r^{\go}(\xi_{i-1},\xi_{i})>0$  for  $i \in \lint 1,|\gG| \rint$.
 For any given ordered pair $(\xi, \xi') \in \gO_{N,k} \times \gO_{N,k}$, we assign a  path $ \gG_{\xi, \xi'}$, whose starting point is $\xi$ and ending point is $\xi'$.

Using \cite[Corollary 13.21]{LPWMCMT}, the spectral gap of the chain can be controlled by a simple quantity  depending on the functional $(\xi, \xi')\mapsto \gG_{\xi,\xi'}$. We say that an unordered pair $e=\{\xi,\xi'\} \subset  \gO_{N,k}$ is an edge if 
$q(e):=\pi^{\go}_{N,k}(\xi) r(\xi,\xi')>0$ (note that  by reversibility $q(e)$ does not depend on the orientation). We write $e\in \gG= (\xi_0,\dots, \xi_{|\gG|})$ if there exists $i\in \lint 1,|\gG|\rint$ such that 
$e=\{ \xi_{i-1},\xi_{i} \} $. We have then (the factor $1/2$ is irrelevant but appears because we are considering unoriented edges rather than oriented one)

\begin{equation}\label{lwbd:gapflow}
 \gap^{\go}_{N,k}\ge  \left(\max_{e} \frac{1}{2q(e)} \sum_{(\xi,\xi')\in \gO_{Nk} \times \gO_{N,k} \ : \ e \in \gG_{\xi,\xi'}}  \pi^{\go}_{N,k}(\xi)\pi^{\go}_{N,k}(\xi') | \gG_{\xi,\xi'} |\right)^{-1}.
\end{equation}
In the proof we describe   a choice for $\gG_{\xi,\xi'}$ which yields a relevant bound for the spectral gap.
Let us  fix  a state $\xi^* \in \Omega_{N,k}$ which has maximal probability, that is such that
\begin{equation}\label{maxenergyconfi}
V^{\go}(\xi^*)= \min_{\xi \in \Omega_{N,k}} V^{\go}(\xi)
\end{equation}
(we make an arbitrary choice if there are several minimizers).
Now to build the path $\gG_{\xi,\xi'}$ we are going to build first a path from $\xi$ to $\xi^*$ and then one from $\xi^*$ to $\xi'$ and then concatenate the two.

\medskip

\noindent We can thus focus on the construction of $\gG_{\xi,\xi^*}$.
Let 
$$m:=d_H(\xi,\xi^*):= \frac{1}{2}\sum_{x=1}^N \left |\xi(x)-\xi^*(x) \right|$$
denote the Hamming distance between $\xi$ and $\xi^*$.
Our first step is to build a sequence $\xi^{(0)},\dots, \xi^{(m)}$ which reduces the Hamming distance in incremental steps that is such that 
\begin{equation}
\begin{cases}
 \xi^{(0)}=\xi \quad \text{ and } \quad \xi^{(m)}=\xi^*,\\
 d_H(\xi^{(i-1)},\xi^{(i)})=1 \quad \quad  \text{ for } i\in \lint 1, m \rint,\\
 d_H(\xi^{(i)},\xi^*)=m-i \quad  \quad \text{ for } i\in \lint 1, m \rint.
\end{cases}
\end{equation}
The choice we make for $\xi^{(0)},\dots, \xi^{(m)}$ is not relevant for the result but let us fix one for the sake of clarity.
Let the sequences $(x_i)_{i=1}^m$ and $(y_i)^m_{i=1}$ be  defined by
\begin{equation}\label{dispreppositions}
\begin{split}
 x_i& \colonequals \min\left\{ x\in \lint 1,N \rint \ : \  \sum_{x=1}^N (\xi(x)-\xi^*(x))_+=i \right\},\\
 y_i& \colonequals \min \left\{ y\in \lint 1,N \rint \ : \  \sum_{y=1}^N (\xi^*(y)-\xi(y))_+=i \right\}.
\end{split}
\end{equation}
These sequences locate the discrepancies between $\xi$ and $\xi^*$.
 Then we define $\xi^{(i)}$ inductively as being  obtained from $\xi^{(i-1)}$ by moving the particle at $x_i$ to $y_i$ which is equivalent to setting
$$ \xi^{(i)}= \xi \wedge \xi^* + \sum_{j=1}^i \ind_{\{y_j\}} +  \sum_{j'=i+1}^m \ind_{\{x_{j'}\}}.$$
Finally, our path from $\xi$ to $\xi^*$ is defined by concatenating paths $\gG^{(i)}$, $i \in \lint 1,m\rint$, linking $\xi^{(i-1)}$ to $\xi^{(i)}$.
We define $\gG^{(i)}= (\xi^{(i)}_0, \dots, \xi^{(i)}_{|x_i-y_i|})$ as   a path of minimal length $|x_i-y_i|$  linking $\xi^{(i)}_0:=\xi^{(i-1)}$ to $ \xi^{(i)}_{|x_i-y_i|}=\xi^{(i)}$. 
To define the intermediate steps, let us assume for notational simplicity (and without loss of generality) that $x_i< y_i$. Moreover,
let $(z_j)_{j=1}^b$  be defined as the  decreasing sequence such that  $$\xi^{(i-1)} \vert_{\lint x_i, y_i\rint}=\ind_{\{z_j\}^b_{j=1}}.$$
We then set $d_j \colonequals  y_i-z_j$ if $j\in \lint 1,b\rint$ and $d_0 \colonequals 0$, and define $(\xi^{(i)}_{\ell})_{\ell=1}^{y_i-x_i}$ by setting if $d_{j-1} < \ell \le d_j$
\begin{equation}\label{moveparticles}
 \xi^{(i)}_{\ell}:=  \xi^{(i-1)}_{\ell}-\ind_{\{z_j\}}+ \ind_{\{z_j+\ell-d_{j-1}\}}.
\end{equation}
In other words, we move the particle at site $z_j$ ($j \ge 1$) to site $z_{j-1}$ (with $z_0=y_i$) starting from $j=1$ until $j=b$. We refer to Figure \ref{fig:flowpath} for a graphical description.

\begin{lemma} \label{lema:paths}
For the path collection $\left(\gG_{\xi,\xi'}\right)$ constructed above, we have  

\begin{equation}
B:=\max_{e} \frac{1}{2q(e)} \sum_{(\xi,\xi')\in \gO_{Nk} \times \gO_{N,k} \ : \ e \in \gG_{\xi,\xi'}}  \pi^{\go}_{N,k}(\xi)\pi^{\go}_{N,k}(\xi') | \gG_{\xi,\xi'} |\le  \ga^{-1} N^2 |\gO_{N,k}|\left(\frac{1-\alpha}{\alpha}\right)^{  N/2}.
\end{equation}
\end{lemma}

 Let us now conclude the proof of Proposition \ref{prop:lwbdgap}.
By \eqref{lwbd:gapflow} and Lemma \ref{lema:paths}, we have
\begin{equation}
\gap_{N,k}^{\go} \ge  \ga N^{-2} |\gO_{N,k}|^{-1}\left(\frac{1-\alpha}{\alpha}\right)^{  -N/2}.
\end{equation}
Observe that 
\begin{equation*}
\max_{\xi, \xi' \in \gO_{N,k}}  \left( V^{\go}(\xi)-V^{\go}(\xi') \right) \le Nk \log \frac{1-\ga}{\ga},
\end{equation*}
and then
\begin{equation}
\min_{\xi \in \gO_{N,k}} \pi_{N,k}^{\go}(\xi) \ge \vert \gO_{N,k} \vert^{-1} \left( \frac{1-\ga}{\ga}  \right)^{-Nk}.
\end{equation}
By \eqref{gapmixrelation}, we have for $\gep \in (0, 1/2)$
\begin{equation}
t_{\Mix}^{N,k, \go}(\gep) \le \ga^{-1} N^{2} \vert \gO_{N,k} \vert \left( \frac{1-\ga}{\ga} \right)^{N}\left( \log \vert \gO_{N,k} \vert +Nk \log \frac{1-\ga}{\ga}  - \log \gep \right).
\end{equation}
\qed

\begin{proof}[Proof of Lemma \ref{lema:paths}]
A first observation is that by construction, our paths are of length smaller than $N^2$.
Let $e$ be an edge and $(\xi,\xi')$ such that $e \in \gG_{\xi,\xi'}$. By symmetry and taking away the factor $1/2$, we can always assume that $e$ belongs to the first part of the path linking $\xi$ to $\xi^*$.
After replacing $| \gG_{\xi,\xi'} |$ by the upper bound and summing over all $\xi'$, we obtain that the quantity we want to bound is exactly 
\begin{equation}\label{darhs}
 \frac{1}{2q(e)} \sum_{(\xi,\xi')\in \gO_{N,k} \times \gO_{N,k} \ : \ e \in \gG_{\xi,\xi'}}  \pi^{\go}_{N,k}(\xi)\pi^{\go}_{N,k}(\xi') | \gG_{\xi,\xi'} |\\  \le N^2 \sum_{\xi\in \gO_{N,k} \ : \ e \in \gG_{\xi,\xi^*}}  \frac{\pi^{\go}_{N,k}(\xi)}{q(e)}.
\end{equation}
Now let $\chi_0(e,\xi)$ denote the first end of $e$ which is visited by the path going from $\xi$ to $\xi^*$. Now simply observing that $q(e)$ is at least  $\ga$ times 
the smallest probability $\pi^{\go}_{N,k}$ of its two end points, we have
\begin{equation}
 \frac{\pi^{\go}_{N,k}(\xi)}{q(e)} \le  \sup_{\xi'\in \gG_{\xi,\xi^*}}  \ga^{-1} e^{V(\xi')-V(\xi)}.
\end{equation}
 Hence using the bound in the sum in \eqref{darhs} we obtain that
\begin{equation}
 \log B \le \log  \ga^{-1} N^2|\gO_{N,k}|+   \suptwo{\xi \in \gO_{N,k}}{\xi'\in \gG_{\xi,\xi^*}} V(\xi')-V(\xi).
\end{equation}
To conclude we only need to prove that for every $\xi \in \gO_{N,k}$ and $\xi'\in \gG_{\xi,\xi^*}$ we have
\begin{equation}\label{zebond}
 V(\xi')-V(\xi) \le \frac{N}{2} \log \left( \frac{1-\alpha}{\alpha}\right).
\end{equation}
This follows simply by inspection from the following observation which follows from our construction and our assumptions.
\begin{itemize}
 \item [(i)] In one step of  $\gG_{\xi,\xi^*}$, $V$ varies at most by $\log \left( \frac{1-\alpha}{\alpha}\right)$ in absolute value.
 \item [(ii)] Along the sequence $(\xi^{(i)})_{i=1}^{m}$,  $V(\xi^{(i)})$ is non-increasing.
 \item [(iii)] Each concatenated path $\gG^{(i)}$ has a length smaller than $N$ (hence each 
$\xi^{(i)}_{\ell}$ is  within   $N/2$ steps of either $\xi^{(i)}$ or $\xi^{(i-1)}$)
so that we have
  $$   \max_{0 \le \ell \le \vert x_i-y_i \vert}( V(\xi^{(i)}_\ell)-V(\xi))\le  \max_{0 \le \ell \le \vert x_i-y_i \vert} \left( V(\xi^{(i)}_\ell)- V(\xi^{(i)})\wedge V(\xi^{(i-1)})\right)\le \frac{N}{2} \log \frac{1-\ga}{\ga}.$$

\end{itemize}

 \begin{figure}[h]
 \centering
   \begin{tikzpicture}[scale=.3,font=\tiny]
     \draw (29,-1) -- (51,-1);     
     \foreach \x in {29,31,...,51} {\draw (\x,-1.3) -- (\x,-1);}
    
     \node[below] at (27,1) {$\xi^{(i-1)}$};

       \node[below] at (29, -1.3) {$x_i$};         
       \draw[fill,red] (29,-1) circle [radius=0.2];
       \draw[fill,black] (31,-1) circle [radius=0.2]; 
        \node[below] at (31, -1.3) {$x_{i+1}$};
         \draw[fill] (41,-1) circle [radius=0.2];
       \node[below] at (41,-1.3) {$x_{i+2}$}; 
       \draw (29.2, -0.5) edge[bend left,->,dashed,red] (46.8, -0.5);

       \node[below,red] at (45, -2){$z_1$};  
       \node[below,red] at (41, -2){$z_2$};        
       \node[below,red] at (39, -2){$z_3$};        
       \node[below,red] at (31, -2){$z_4$};
       \node[below,red] at (29, -2){$z_5$};

      \draw (45.2, -0.8) edge[bend left,->] (46.8, -0.8);             
      \node[above] at (46,-0.8) {$1$};
      \draw[fill] (39,-1) circle [radius=0.2];
      
      \draw (41.2, -0.8) edge[bend left,->] (42.8, -0.8);             
       \node[above] at (42,-0.8) {$2$}; 
      \draw (43.2, -0.8) edge[bend left,->] (44.8, -0.8);             
       \node[above] at (44,-0.8) {$3$}; 
      \draw (39.2, -0.8) edge[bend left,->] (40.8, -0.8);             
       \node[above] at (40,-0.8) {$4$}; 
      \draw (31.2, -0.8) edge[bend left,->] (32.8, -0.8);             
       \node[above] at (32,-0.8) {$5$}; 
       \draw (33.2, -0.8) edge[bend left,->] (34.8, -0.8);             
       \node[above] at (34,-0.8) {$6$}; 
        \draw (35.2, -0.8) edge[bend left,->] (36.8, -0.8);             
       \node[above] at (36,-0.8) {$7$}; 
        \draw (37.2, -0.8) edge[bend left,->] (38.8, -0.8);             
       \node[above] at (38,-0.8) {$8$};  
       \draw (29.2, -0.8) edge[bend left,->] (30.8, -0.8);             
       \node[above] at (30,-0.8) {$9$};   
       
       \draw[fill] (45,-1) circle [radius=0.2];        
       \draw[fill] (51,-1) circle [radius=0.2];

     \draw (29,-6) -- (51,-6);
     \foreach \x in {29, 31,...,51} {\draw (\x,-6.3) -- (\x,-6);}

     \node[below] at (27,1-5) {$\xi^{(i)}$};

    \draw[fill] (51,-1-5) circle [radius=0.2]; 
   
    \draw[fill,blue] (47,-1-5) circle [radius=0.2]; 
    \draw[fill] (45,-1-5) circle [radius=0.2];
    
    \node[below] at (47, -1.3-5) {$y_i$};
                          
    \draw[fill] (39,-1-5) circle [radius=0.2];  
    \draw[fill,black] (31,-1-5) circle [radius=0.2]; 
          \draw[fill] (41,-1-5) circle [radius=0.2];

    \node[below] at (37, -8) {(L)};

    \node[below] at (63, -8) {(R)};                        

     \draw [->](55,1) -- (65,1)node[right] {$\ell$};
     \draw [->] (55,-7)--(55,2) node[right] {$V(\xi^{(i)}_{\ell})$};
     
     \draw[fill] (55,-2) circle [radius=0.1];  
     \draw[fill] (56,-3.2) circle [radius=0.1];  
     \draw[fill] (57,-4) circle [radius=0.1];  
     \draw[fill] (58,-5) circle [radius=0.1];  
     \draw[fill] (59,-4.3) circle [radius=0.1];  
     \draw[fill] (60,-3) circle [radius=0.1];  
     \draw[fill] (61,-4) circle [radius=0.1];  
     \draw[fill] (62,-5) circle [radius=0.1];  
     \draw[fill] (63,-6) circle [radius=0.1]; 
     \draw[fill] (64,-5) circle [radius=0.1];  
     \draw[dashed](55,-2)--(56,-3.2)--(57,-4)--(58,-5)--(59,-4.3)--
     (60,-3)--(61,-4)--(62,-5)--(63,-6)--(64,-5);

 \end{tikzpicture} \label{fig:flowpath} 
   \caption{ 
    A bold circle represents a particle, and a particle at the same site for the configurations $\xi^{(i-1)}$ and $\xi^{(i)}$ is colored black. Otherwise, it is  red or blue.
   (L) A graphical description of the movements of the particle at site $x_i$  of $\xi^{(i-1)}$ to the empty site $y_i$ and the numbers above the arrows are the relative order of the movements.
   (R) We draw the graph of $(\ell, V(\xi_{\ell}^{(i)}))_{\ell}$. } 
 \end{figure}
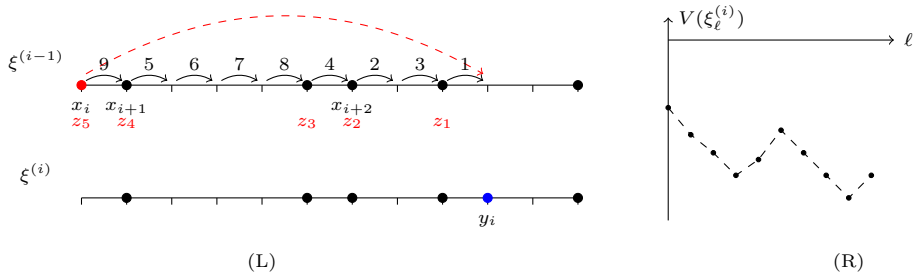

 \end{proof}

\section{Lower bounds on the mixing time}\label{sec:lowbdmix}

 Theorem \ref{th:mixinglb} contains three separate lower bounds. The first one is a consequence of Proposition \ref{prop:ulb}. 
 In this section, we are going to prove the two remaining bounds which are restated below as Propositions \ref{prop:oneparticle} and \ref{prop:flow} respectively. The proof of these propositions rely on the  two mechanisms exposed in Subsection \ref{subsec:interpconject}: The potential barrier created by rare fluctuations of $V^{\go}$ (cf. Proposition \ref{lema:createtrap}) has the effect of trapping individual particles and slowing down the particle flow.

\subsection{A lower bound from the position of the first particle}

\begin{proposition}\label{prop:oneparticle}
We have with high probability
\begin{equation}\label{consequence}
 t_{\Mix}^{N,k, \go} \geq  [N(\log N)^{-2}]^{\frac{1}{\gl}}
\end{equation}

\end{proposition}

\begin{proof}
We let $y_1(\go)> x_1(\go)$ be such that 
$V(y_1(\go))-V(x_1(\go))$ is maximized within $1\le x\le y \le N/4$ (the event that $V^{\go}$ is non-increasing on $\lint 1, N/4\rint$ is unlikely, and then it can be ignored). 
We are going to prove that w.h.p.\ 
\begin{equation}\label{cause}
t_{\Mix}^{N,k, \go} \geq \frac 1 {2e} e^{V(y_1)-V(x_1)}-1.
\end{equation}
As a consequence 
of Proposition \ref{lema:createtrap} (applied to the segment $\lint 1, N/4\rint$), we have w.h.p.\ 
$$ V(y_1)-V(x_1) \ge \frac{1}{\gl} \log N-\frac{2}{\gl} \log \log N+\log 20,$$
so that \eqref{consequence} follows from \eqref{cause}. Recall the notation \eqref{partiposition}. Then considering that  $\xi_{\min}$ should be the worst initial condition, we observe that 
\begin{equation}
d_{N,k}^{\go}(t) \ge\bP \left[ \bar \gs^{\min}_t(1) \le \frac{N}{4} \right]-\pi_{N,k}^{\go} \left(\bar \xi(1)\le N/4\right). 
\end{equation}
As a consequence of  Lemma \ref{lema:equilleftright} we have w.h.p. $$\pi_{N,k}^{\go}(\bar \xi(1)\le N/4)\le 1/4.$$
To have an estimate on the mixing time, we must prove that $\bP \left[ \bar \gs^{\min}_t(1) \le \frac{N}{4} \right]\ge 1/2$. We define
$$\tau_{y_1} \colonequals  \inf \left\{ t \ge 0: \bar \gs^{\min}_t(1)=y_1  \right\}.$$
We have 
\begin{equation}
 \bP \left[ \bar \gs^{\min}_t(1) > \frac{N}{4} \right]\le \bP[\tau_{y_1}\le t].
\end{equation}
We are going to show that 
\begin{equation}
  \bP[\tau_{y_1}\le t] \le e (t+1) e^{V(x_1)-V(y_1)},
\end{equation}
which is sufficient for us to conclude that \eqref{cause} holds. 
Using the graphical construction (with an enlargement of the probability space to sample the initial condition) we can couple $\gs^{\min}_t$ with $X_t^{\pi}$ a random walk on the interval $\lint 1, y_1 \rint$ with transitions rates given by $q^{\go}_{y_1}$  (cf. \eqref{qngo}) and starting with an initial distribution sampled from the equilibrium measure $\pi_{y_1, 1}^{\go}$, in such a way that 
$$\forall t\le \tau_{y_1},  \quad   \bar \gs^{\min}_t(1) \le X_t^{\pi}.$$
Setting $\tilde \tau_{y_1} \colonequals \inf \left\{t \ge 0: X_t^{\pi} =y_1 \right\}$,
we then  have 
\begin{equation}\label{onepartbound}
  \bP \left[ \tau_{y_1} \le t \right]\le \bP \left[ \tilde \tau_{y_1} \le t \right].
\end{equation} 
We define the occupation time 
$$u(t) \colonequals \int_0^t \ind_{\{ y_1\}}(X_s^{\pi}) \dd s. $$
We have
 \begin{equation}\label{lwbd:hittime}
\bE \left[ u(t+1) \right]\ge \bP[ u(t+1)\ge 1] \ge 
\bP[ \tilde \tau_{y_1}\le t] \bP \left[ \forall s \in [0,1]: X^{\pi}_{\tilde \tau_{y_1}+s}=y_1  \right] \ge e^{-1} \bP[ \tilde \tau_1\le t],
\end{equation}
where in the last inequality we use the strong Markov property.
 As  the process $(X_t^{\pi})_{t \ge 0}$ is stationary,  
 $$\bE[ u(t+1)]= (t+1)\pi^{\go}_{y_1,1}(y_1)\le  (t+1) e^{V(y_1)-V(x_1)},$$
 which allows to conclude that 
 \begin{equation}
  \bP[ \tilde \tau_1\le t]\le e(t+1)e^{V(y_1)-V(x_1)}.
 \end{equation}

\end{proof}

\subsection{A lower bound derived from flow consideration}

 Let us now derive the third bound which is necessary to complete the proof of Theorem \ref{th:mixinglb}.

\begin{proposition}\label{prop:flow}
There exists  a positive  constant $c=c(\ga, \bbP)$   such that w.h.p.\ we have
\begin{equation}\label{lowbdmix}
 t_{\Mix}^{N,k}\geq  c  k  N^{\frac{1}{2\gl}} (\log N)^{-2\left(1+\frac{1}{\gl}\right)} .
\end{equation}
\end{proposition}

To prove the above result, 
we adopt the strategy developed in \cite[Proposition 4.2]{schmid2019mixing} by investigating the flow of particles through a \textit{slow} segment of size of order $(\log N)$ where the drift of the random environment points to the left.
This flow of particles is controlled via a comparison with a  boundary driven exclusion process.

\medskip

In \cite{schmid2019mixing} the \textit{slow} segment is selected to be such that 
$\go_x<1/2$ for \textit{every site}. It has the advantage of simplifying the computation since it allows for comparison with the homogeneous exclusion process for which computation has been performed in \cite{blythe2000exact}.
Our approach brings an improvement by selecting the slow segment based on the potential function $V^{\go}$. The relevant quantity that limits the flow is the worst potential barrier that the particles have to overcome. Proposition \ref{lema:createtrap} allows to identify the worst potential barrier in the system.
We let $x_2(\go)\le y_2(\go)$ be the smallest elements of $\lint N/2, 3N/4\rint$ such that 
$$V^{\go}(y_2)-V^{\go}(x_2)=\max_{N/2\le x\le y\le 3N/4} \left( V^{\go}(y)-V^{\go}(x) \right).$$
According to Proposition \ref{lema:createtrap} we have w.h.p.\
\begin{equation}\label{trapenvironm}
 V(y_2)-V(x_2)\ge \frac{1}{\gl} \left(\log N- 2\log \log N \right)   \text{ and } y_2-x_2\le q_N. 
\end{equation}

 In order to illustrate how the mixing time can be controlled using the flow of particles, 
 we start with a simple lemma.
Let $J_t$ denote the number of particles on the last portion of the segment,
\begin{equation}\label{def:Jt}
 J_t:=  \sum_{x\ge y_2+1} \sigma^{\min}_t(x).
\end{equation}

\begin{lemma}\label{lema:distflow}
 For any $\gep>0$, we have with high probability for every $t\ge 0$.
\begin{equation}
 d_{N,k}^{\go}(t)\ge 1- \frac{4\bE[J_t]}{k}-\gep.
\end{equation}
\end{lemma}

\begin{proof}
Setting 
$\cB:= \left\{ \xi \in \gO_{N,k} \ : \ \sum_{x\ge y_2+1} \xi(x)< k/4 \right\}$, 
 we have
 \begin{equation}
   d_{N,k}^{\go}(t) \ge   \| P^{\xi_{\min}}_t - \pi_{N,k}^{\go} \|_{\TV}
   \ge \bP[ \sigma^{\min}_t \in \cB]-  \pi_{N,k}^{\go}(\cB).
 \end{equation}
By Lemma \ref{lema:equilleftright}, the second term  is smaller than $\gep$ with high probability.
Concerning the first term, 
we have by Markov's inequality
\begin{equation}
  \bP \left[ \sigma^{\min}_t \in \cB \right]= 1- \bP \left[ J_t\ge k/4 \right] \ge 1- \frac{4\bE[J_t]}{k}.
\end{equation}

\end{proof}

\begin{figure}[h]
 \centering
   \begin{tikzpicture}[scale=.4,font=\tiny]
     \draw (25,-1) -- (53,-1);
     
     \foreach \x in {25, 27,...,51,53} {\draw (\x,-1.3) -- (\x,-1);}
     \draw[fill] (51,-1) circle [radius=0.2];   
     \node[below] at (25,-1.3) {$x_2$};
     \node[below] at (51,-1.3) {$y_2$};
     \node[below] at (53.5,-1.3) {$y_2+1$};

       \draw (23, 1) edge[bend left,->] (25, -0.8);
       \node[below] at (24.7, 1.6) {$\go_{x_2-1}$};

       \draw[fill] (23,1) circle [radius=0.2];  
       \node[below] at (29, -1.3) {$x$};
       \node[below] at (28, 0.5) {$1-\go_x$};
       \node[below,red] at (30, 0.5) {$\times$};    

       \draw (51, -0.7) edge[bend left,->] (52.7, 1);
       \draw (51, -0.7) edge[bend right,->] (49, -0.7);
       \node[below] at (51.5,1.8) {$\go_{y_2}$};
       \node[below] at (49.5,0.8) {$1-\go_{y_2}$};
       \draw (52.7,2)--(52.7,-1)--(53.3,-1)--(53.3,2);
        \foreach \x in {-1,-0.55,...,1} {\draw[fill] (53,\x) circle [radius=0.2]; }
       
       \draw[fill] (29,-1) circle [radius=0.2];
   
       \draw (29, -0.8) edge[bend right,->] (27.2, -0.8);
       \draw (29, -0.8) edge[bend left,dashed,->, red] (30.8, -0.8);
       \draw[fill] (31,-1) circle [radius=0.2]; 
       
       \draw (39, -0.8) edge[bend left,->] (40.8, -0.8);   
       \draw (39, -0.8) edge[bend right,->] (37.2, -0.8);      
       \node[below] at (39,-1.3) {$y$};
       \node[below] at (40,0.5) {$\go_y$};
       \node[below] at (38,0.5) {$1-\go_y$};

       \draw[fill] (39,-1) circle [radius=0.2];
         
       \draw (45, -0.8) edge[bend right, ->] (43.2, -0.8);  
       \draw (45, -0.8) edge[bend left,->,dashed,red] (46.8, -0.8);      
       \node[below, red] at (46,0.5) {$\times$};
       \node[below] at (44,0.5) {$1-\go_z$};
       \node[below] at (45,-1.3) {$z$};
   
       \draw[fill] (45,-1) circle [radius=0.2];   
       \draw[fill] (47,-1) circle [radius=0.2];       
 \end{tikzpicture}
   \caption{ A graphical representation of the boundary driven process: a bold circle represents a particle, and the number above every arrow represents the jump rate while a red $"\times"$ represents an inadmissible jump. In addition, the site $y_2+1$ can accommodate infinite many particles and all particles at site $y_2+1$ stay put. } \label{fig:bounddriven} 
 \end{figure}
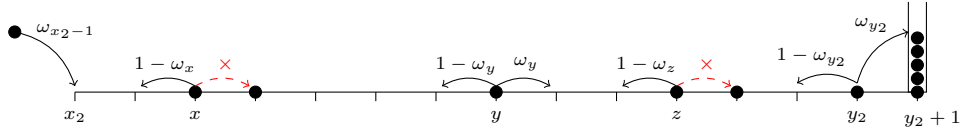

Now we can control $\bE[J_t]$ by comparing our system with one in which the particles flow faster. 
We consider a process on a different state space
\begin{equation}
 \tilde \gO_{x_2,y_2} \colonequals  \left\{  \xi : \lint x_2, y_2+1 \rint \to \bbZ_+ \ : \ 
 \forall x\in \lint x_2, y_2 \rint, \    \xi(x)\in \{ 0,1\} \right\}.
\end{equation}
Under this new process the particles follow the exclusion dynamics in the bulk but new rules are added at the boundaries. If $\xi(x_2)=0$ then a particle is added at site $x_2$ with rate one. At the other end of the segment  particles can jump from site $y_2$ to site $y_2+1$ without respecting the exclusion rule (\textit{i.e.}, the site $y_2+1$ is allowed to contain arbitrarily many particles) and particles at site $y_2+1$ remain there forever.
We define the generator of the process to be (for $f: \tilde \gO_{x_2,y_2} \mapsto \bbR $)
\begin{multline}\label{generator:bddirven}
\tilde {\mathfrak{L}}_{x_2,y_2}^{\omega}f(\xi) \colonequals  \sum_{z=x_2}^{y_2-1}
 r^{\go}(\xi,\xi^{z,z+1})\big[f(\xi^{z,z+1})-f(\xi)\big]\\
+ \go_{x_2-1} \ind_{\{\xi(x_2)=0\}}\big[f(\xi+ \delta_{x_2})-f(\xi)\big]+ \go_{y_2}\ind_{\{\xi(y_2)=1\}} \big[f(\xi-\delta_{y_2}+\delta_{y_2+1})-f(\xi)\big],
\end{multline}
where $r^{\go}$ is defined in \eqref{jumprateSEP}. We refer to Figure \ref{fig:bounddriven} for a graphical description. 
We  let $(\tilde \sigma^{\xi}_t)_{t \ge 0}$ denote the corresponding process starting from an initial condition $\xi \in \tilde \gO_{x_2,y_2}$.

\begin{lemma}\label{lema:modifiedflow}
 Let $\bf 0$ denote the configuration with all sites in $\lint x_2, y_2+1 \rint$ being empty, and then  we have 
 \begin{equation}
    J_t \le  \tilde \sigma^{\bf 0}_t(y_2+1),
 \end{equation}
 where $J_t$ is defined in \eqref{def:Jt}.
\end{lemma}

\begin{proof}
The process $(\tilde \sigma^{\bf 0}_t)_{t \ge 0}$ can be constructed together with $(\sigma^{\min}_t)_{t \ge 0}$ on the same probability space using the graphical construction of Section \ref{subsec:cancoupling} (with the obvious adaptation of the construction to fit the boundary condition for $(\tilde \sigma^{\bf 0}_t)_{t \ge 0}$ using the same clocks $(T^{(x)}_n)_{x,n \in \bbN}$ and auxiliary variables $(U^{(x)}_n)_{x,n \in \bbN}$ for the two processes.
It can then be checked by inspection that for every $t\ge 0$
\begin{equation}\label{wooz}
 \forall x\in \lint x_2, y_2+1\rint \quad  \sum_{z=x}^{N} \sigma^{\min}_t(z)\le  \sum_{z=x}^{y_2+1} \tilde \sigma^{\bf 0}_t(z).
\end{equation}
Since the above inequality is satisfied at $t=0$, it is sufficient to check that it is conserved by any update of the two processes.
The result then just corresponds to the case $x=y_2+1$.
\end{proof}

\begin{proposition} \label{prop:partijumpout}
There  exists a  constant $C=C(\ga,\bbP)$ such that for all $t \ge 0$ w.h.p. we have
\begin{equation}
\bE \left[ \tilde \sigma_t(y_2+1)\right]\leq t C   N^{-\frac{1}{2\gl}} (\log N)^{2(1+\frac{1}{\gl})}.
\end{equation}
\end{proposition}

\noindent With Proposition \ref{prop:partijumpout} whose proof is detailed in the next subsection, we are ready for the proof of Proposition \ref{prop:flow}.

\begin{proof}[Proof of Proposition \ref{prop:flow}] By Lemma \ref{lema:distflow} and Lemma \ref{lema:modifiedflow}, we have
\begin{equation}\label{dist:lwbdflow}
 d_{N,k}^{\go}(t)\ge \frac{7}{8}- 4\frac{\bE[\tilde \sigma_t(y_2+1)]}{k}.
\end{equation}
By Proposition \ref{prop:partijumpout}, we take $$t=\frac{1}{8C} k N^{\frac{1}{2\gl}} (\log N)^{-2(1+\frac{1}{\gl})}$$
in \eqref{dist:lwbdflow} to conclude the proof.

\end{proof}

\subsection{Proof of Proposition \ref{prop:partijumpout}.}\label{subsec:bounddriven}

 Note that $\tilde \sigma^{\bf 0}_t(y_2+1)$ is  a  superadditive ergodic sequence. 
To see this we let $\vartheta_s$ denote the time shift operator on the graphical construction variables. Recalling \eqref{defizero} we set
\begin{equation}
  \vartheta_s(  (T^{(x)}_i, U^{(x)}_i)_{x\in \bbZ, i\ge 1}):=   \left( T^{(x)}_{i+i_0(x,s)}-s, U^{(x)}_{i+i_0(x,s)}\right)_{x\in \bbZ, i\ge 1}.
\end{equation}
Now we observe that the graphical construction preserves the order $\preccurlyeq$ on $\tilde\gO_{x_2,y_2}$ defined by 
\begin{equation}
 \xi \preccurlyeq \xi'  \quad  \text{ if } \quad   \forall x\ge x_2, \quad  \sum_{z=x}^{y_2+1} \xi(z)\le  \sum_{z=x}^{y_2+1} \xi'(z).
\end{equation}
 Hence comparing the dynamic in the interval $[s,s+t]$ with that starting from ${\bf 0}$ at time $s$, we obtain
\begin{equation}\label{subbaa}
 \tilde \sigma^{\bf 0}_{s+t}(y_2+1)\ge  \tilde \sigma^{\bf 0}_{s}(y_2+1)+ (\vartheta_s \circ \tilde \sigma)_t^{\bf 0}(y_2+1). 
 \end{equation}
  Since the shift operator $\vartheta_s$ on $(T,U)$ is ergodic, we obtain from Kingman's subbadditive ergodic Theorem \cite{kingman73} (continuous time version)  that
\begin{equation}\label{ergodic}
 \bE \left[\tilde \sigma^{\bf 0}_{t}(y_2+1)\right] \le t  \left[\lim_{s \to \infty}\frac{1}{s}  \tilde \sigma^{\bf 0}_{s}(y_2+1)\right].
\end{equation}
Letting $\cN_{s}:=\sum_{x=x_2}^{y_2} \tilde \sigma^{\bf 0}_s(x)$ denote the number of  mobile particles in the system (particles at site $y_2+1$ which have stopped moving are not counted), we have
\begin{equation}
  \tilde \sigma^{\bf 0}_t(y_2+1)= \sum_{s\in(0,t]} \ind_{\{\cN_{s}<\cN_{s_-}\}}.
\end{equation}
Letting $(\cT_n)_{n\ge 1}$ denote the sequence of time at which $\cN_{t}<\cN_{t_-}$ (in increasing order), we have 
\begin{equation}\label{npartjumpout}
 \lim_{s \to \infty}\frac{1}{s}  \tilde \sigma^{\bf 0}_{s}(y_2+1)= \lim_{n\to \infty} \frac{n}{\cT_n}.
\end{equation}
 Similarly to \eqref{subbaa}, using preservation of order and the fact that 
$$ \sigma^{\bf 0}_{s} \preccurlyeq  [ \sigma^{\bf 0}_{s}(y_2+1)+y_2-x_2+1]\ind_{\{y_2+1\}},$$
we have for every $s>0$
\begin{equation}\label{subaddz}
 \tilde \sigma^{\bf 0}_{s+t}(y_2+1)
 \le \tilde \sigma^{\bf 0}_{s}(y_2+1)+ (\vartheta_s \circ \tilde \sigma)_t^{\bf 0}(y_2+1)+y_2-x_2+1.
\end{equation}
Now as a consequence of  \eqref{subaddz}  
\begin{equation}
 \cT_{l+y_2-x_2+2}\ge  \cT_l + \vartheta_{\cT_l}\circ  \cT_{1}.
\end{equation}
Since $\cT_l$ is a stopping time with respect to $(\cF_t)_{t \ge 0}$ (recall \eqref{filtrations}), by the strong Markov property
$\vartheta_{\cT_l}\circ  \cT_{1}$ is independent of $\cT_1$ and has the same distribution.
Iterating the process we obtain that 
\begin{equation}
  \cT_{(r-1)(y_2-x_2+2)+1}\ge \cT^{(1)}_1+\cdots+  \cT^{(r)}_1
\end{equation}
where $(\cT^{(a)}_1)_{a=1}^r$ is a sequence of IID copies of $ \cT_{1}$. 
This yields that 
\begin{equation}\label{onejumpout}
 \liminf_{n\to \infty} \frac{\cT_n}{n}\ge  \frac{1}{y_2-x_2+2} \bE\left[\cT_{1} \right].
\end{equation}
Finally let us compare $(\tilde \sigma^{\bf 0}_t)_{t \ge 0}$ with  $(\tilde \sigma^{'}_t)_{t \ge 0}$  starting from 
another initial condition. Now we specify the initial condition.  Let us first choose the number of particle by setting
\begin{equation}\label{choosek}
\begin{split}
\gL(\go)& \colonequals \left\{ x\in \lint x_2,y_2 \rint \ : \ V(x)  \le  [V(y_2)+V(x_2)]/2 \right\},\\
  k'(\go) &   \colonequals \# \gL(\go).
  \end{split}
\end{equation}
We let $(\tilde \sigma^{'}_t)_{t \ge 0}$ be the dynamic with generator \eqref{generator:bddirven} and initial configuration $\tilde \sigma'_0$ is obtained by setting $ \tilde \sigma'_0(y_2+1)=0$ and
sampling $\pi_{[ x_2,y_2],k'}^{\go}$ (the invariant probability measure for the exclusion process on the segment $\lint x_2, y_2\rint$ with $k'$ particles) to set the values of $(\tilde \sigma'_0(x))_{x\in \lint x_2,y_2\rint}$.
Using monotonicity again we have
\begin{equation}\label{thacompa}
\cT_{1} \ge \inf\{ t \ge 0 :  \tilde \sigma'_t(y_2+1)=1 \}
\ge \inf\{ t \ge 0 :  \tilde \sigma'_t(x_2)=0 \ \text{ or } \   \tilde \sigma'_t(y_2)=1 \} \equalscolon \cT'.
\end{equation}
Now let us observe that until time $\cT'$, the process $(\tilde \sigma'_t)_{t \ge 0}$ (or rather, its restriction to $\lint x_2,y_2\rint$) coincides with the exclusion process on the segment $\lint x_2, y_2\rint$ with $k'$ particles. Using this we can prove the following (the proof is postponed to the end of the section).

\begin{lemma} \label{lema:hittingtime}
We have
\begin{equation}
\bE \left[\cT' \right]  \ge  \frac{1}{16e^2 (y_2-x_2)} e^{\frac{V(y_2)-V(x_2)}{2}}.
\end{equation}
\end{lemma}
Let us now conclude the proof of Proposition \ref{prop:partijumpout}.
Combing \eqref{ergodic}, \eqref{npartjumpout}, \eqref{onejumpout} and \eqref{thacompa}, we have
\begin{equation}
\bE \left[ \tilde \gs^{\bf 0}_t(y_2+1)\right] \le t   \left[\lim_{s \to \infty}\frac{1}{s}  \tilde \sigma^{\bf 0}_{s}(y_2+1)\right] \le  \frac{ t (y_2-x_2+2)}{\bE \left[ \cT_1\right]}\le 
\frac{ t (y_2-x_2+2)}{\bE \left[ \cT'\right]}.
\end{equation}
Using  Lemma \ref{lema:hittingtime},
we obtain
\begin{equation}
\bE \left[ \tilde \gs^{\bf 0}_t(y_2+1)\right] \le t 16e^2 (y_2-x_2+2)^2 e^{-\frac{V(y_2)-V(x_2)}{2}}.
\end{equation}
By \eqref{trapenvironm}, we have  w.h.p. \
\begin{equation}
\bE \left[ \tilde \gs^{\bf 0}_t(y_2+1)\right] \le t 16e^2 (q_N+2)^2 N^{-\frac{1}{2 \gl}} (\log N)^{\frac{1.}{\gl}}.
\end{equation}

\qed

\begin{proof}{Proof of Lemma \ref{lema:hittingtime}}
With a small abuse of notation, in this proof $(\tilde \gs'_t)_{t \ge 0}$ denotes the exclusion process on the segment $\lint x_2,y_2\rint$ with $k'$ particles starting from stationarity.
Since $\bE \left[ \cT' \right] \ge t \bP\left[ \cT' > t  \right]$, our goal is to provide a lower bound on $\bP\left[ \cT' > t  \right]$. 
We  define
\begin{equation}
\begin{gathered}
\cB_1 \colonequals \left\{ \xi \in  \gO_{\lint x_2, y_2 \rint, k'}: \xi(x_2)=0  \right\},  \\
 \cB_2 \colonequals \left\{\xi \in  \gO_{\lint x_2, y_2 \rint, k'}: \xi(y_2)=1 \right\}.
\end{gathered}
\end{equation}
Using the strong Markov property at $\cT'$ and the fact that jumping rates for particles  are bounded from above by one at every site, we have 
\begin{equation*}
\bP \left[\forall t \in [\cT', \cT'+1], \ \tilde \gs_t' \in \cB_1 \cup \cB_2   \right]\ge e^{-2}. 
\end{equation*}
Using independence as  in \eqref{lwbd:hittime}, we have
\begin{equation}\label{hittinghalved}
\bP \left[ \cT' \le  t \right] \le e^2 (t+1) \pi_{[ x_2,y_2],k'}^{\go} \left(\cB_1 \cup \cB_2 \right).
\end{equation}
We now head to provide an upper bound on $\pi_{[ x_2,y_2], k'}^{\go}(\cB_1)$. 
Recalling the definition of $\gL$ in \eqref{choosek}, we observe that when $ \xi \in \cB_1$, since $x_2\in \gL$ and there are $k'$ particles, there must be a particle in $\gL^{\complement} \colonequals \lint x_2,y_2 \rint \setminus \gL$.
Let $\mathtt{R}(\xi)$ be the position of the rightmost such particle 
\begin{equation*}
\mathtt{R}(\xi) \colonequals \sup  \left\{ z \in \gL^{\complement}:\  \xi(z)=1 \right\},
\end{equation*}
and set for $z\in \gL^{\complement}$
\begin{equation*}
\cB_{1,z} \colonequals \left\{ \xi \in \cB_1: \mathtt{R}(\xi)=z \right\}.
\end{equation*}
By moving the particle from site $z$ to site $x_2$ as in \eqref{exchange01}, 
we obtain
\begin{equation*}
\pi_{[ x_2,y_2], k'}^{\go}\left( \cB_{1,z} \right)= \sum_{\xi \in  \cB_{1,z}  }\pi_{[ x_2,y_2], k'}^{\go}(\xi^{x_2,z})e^{-V(z)+V(x_2)} \le e^{-V(z)+V(x_2)} \le e^{-\frac{V(y_2)-V(x_2)}{2}},
\end{equation*}
and then
\begin{equation}\label{trapprob:left}
\pi_{[ x_2,y_2], k'}^{\go}(\cB_1)= \sum_{z \in  \gL^{\complement}} \pi_{[ x_2,y_2], k'}^{\go}(\cB_{1,z}) \le (y_2-x_2) e^{-\frac{V(y_2)-V(x_2)}{2}}.
\end{equation}
Similarly, we can obtain
\begin{equation}\label{trapprob:right}
\pi_{[ x_2,y_2 ], k'}^{\go}(\cB_2) \le (y_2-x_2) e^{-\frac{V(y_2)-V(x_2)}{2}}.
\end{equation}
Combining \eqref{trapprob:left} with \eqref{trapprob:right}, in \eqref{hittinghalved} we take $$t= \frac{1}{4e^2(y_2-x_2)}e^{\frac{V(y_2)-V(x_2)}{2}}-1$$
to obtain
\begin{equation}
\bE \left [\cT' \right] \ge  \frac{1}{2} \left( \frac{1}{4e^2(y_2-x_2)}e^{\frac{V(y_2)-V(x_2)}{2}}-1 \right) \ge \frac{1}{16e^2(y_2-x_2)}e^{\frac{V(y_2)-V(x_2)}{2}}.
\end{equation}

\end{proof}

\section{Upper bound on the mixing time}\label{sec:upbdmix}
 This section is dedicated to the proof of Theorem \ref{th:mixingub}.
First in Section \ref{redhit} we are going to reduce the problem to the estimation of the hitting time of $\xi_{\max}$.
Afterwards using Proposition \ref{PWcoro} and Proposition \ref{prop:censor} we are going to provide estimate of this hitting time using a modified censored dynamics. 
First in Section \ref{subsec:smallk} we treat the case of $k\le q_N$ which is a bit simpler and treat the more general case $q_N<k\le N/2$ in Section \ref{morethank}.

\subsection{Deducing the mixing time from the hitting time of the maximal configuration}\label{redhit}
 Let us first show that the study of the mixing time can be reduced to that of the probability of hitting the configuration $\xi_{\max}$ starting from the other extremal configuration $\xi_{\min}.$

\begin{proposition}\label{prop:reduhit}
 We have for every $t>0$ and $n\in \bbN$
 
 \begin{equation}
  d_{N,k}^{\omega}(nt)\le   \left(1- P_t(\xi_{\min},\xi_{\max})\right)^n.
 \end{equation}

\end{proposition}

\begin{proof}

We have (see for instance \cite[Lemma 4.10]{LPWMCMT})
\begin{equation}\label{distaz}
d_{N,k}^{\omega}(t)\le \bar d_{N,k}^{\omega}(t):=\max_{\xi,\xi'}\| P^{\xi}_t-P^{\xi'}_t \|_{\TV}
\le \max_{\xi,\xi'} \bP\left[ \sigma^{\xi}_t\ne\sigma^{\xi'}_t  \right]
\end{equation}
Using  the monotonicity under the graphical construction (cf. Proposition \ref{monocouple}) 
 for all $\xi \in \Omega_{N,k}$ and $t\ge 0$ we have 
 \begin{equation*}
 \sigma_{t}^{\min}\le \sigma_{t}^{\xi} \le \sigma^{\max}_t,
 \end{equation*}
where $\sigma^{\min}$ and $\sigma^{\max}$ are starting from the extremal conditions $\xi_{\min}$ and $\xi_{\max}$ in \eqref{minmaxconf}.
As a consequence for arbitrary $\xi$ and $\xi'$ with $\tau' \colonequals \inf\{t \ge 0 \ : \sigma_{t}^{\xi}= \sigma^{\xi'}_t\}$, we have 
\begin{equation}
 \forall t \ge \tau', \quad \sigma^{\xi}_t=\sigma^{\xi'}_t.
\end{equation}
On the other hand we have
 \begin{equation}
 \tau' \ge \tau \colonequals \inf \left\{t\geq 0: \sigma_t^{\min} =\xi_{\max}\right\}.
 \end{equation}
 Therefore \eqref{distaz} implies that 
 \begin{equation}
 d_{N,k}^{\omega}(t)\leq \mathbf{P} \big(\tau >t\big).
 \end{equation}
 Using again the 
 Markov property and the monotonicity in Proposition \ref{monocouple}, we have for any positive integer $n$
 \begin{equation}\label{iteratehitting}
 \bP \left(\tau >n t \right)\leq 
 \bP \left(\sigma_{it}^{\min}\neq \xi_{\max}, \forall i \in \lint 1, n \rint \right)
 \leq  \bP \left(\sigma_{t}^{\min}\neq \xi_{\max} \right)^n.
 \end{equation}
 
 \end{proof}

\subsection{The case $k_N\le q_N$}\label{subsec:smallk}

 Before stating the main result of this section, let us present a strategy to bound $P_t(\xi_{\min},\xi_{\max})$ from below. We present in the process a few key technical lemmas whose proof is presented  in the next subsection.
We consider environment within the following event
\begin{equation}
  \cA_N \colonequals \left \{\go:  \maxtwo{ 1\le x \le y \le N}{y-x\ge q_N} \left( V(y)-V(x) \right) \le -3\log N \right\}.
\end{equation}
Note that by Proposition \ref{lema:createtrap}, this is an high probability event.
    The event $\cA_N$ ensures that on segments of length $4q_N$, at equilibrium the particles concentrate on the right half of the segment with high probability. It also ensures that with high probability the last site is occupied by a particle.

\begin{lemma}\label{lema:LhalfRone}
If $\go\in \cA_N$, then we have for any $x\in \lint0,N-4q_N\rint$ and any $k\le q_N$,
 \begin{equation}
 \begin{split}
  \pi_{[x+1,x+4q_N],k}^{\go}\left[\bar \xi(1)\le x+2q_N\right]\le  2q_N^2N^{-3},\\
\pi_{[x+1,x+4q_N], q_N}^{\go}\left[\xi(x+4q_N)=0\right]\le 3q_N N^{-3}.
\end{split}
\end{equation}

\end{lemma}

\noindent Our second technical lemma is  a direct consequence of  Proposition \ref{prop:lwbdgap}.
It allows to bound the mixing time of the system for each of the intervals of length $4 q_N$ in a quantitative way.
We define 
$$T=T_N \colonequals 80 \ga^{-1} q_N^{4} \binom{4q_N}{q_N}\left(\frac{1-\alpha}{\alpha}\right)^{2q_N} \log \left( \frac{1-\ga}{\ga} \right).$$
 The following result is obtained by  taking $\gep=N^{-3}$ in Proposition \ref{prop:lwbdgap}.
\begin{lemma}\label{lema:mixN3}
Under the assumption \eqref{detele}
we have for all $k \le q_N$, all $\go$ and all $x \in \lint 0, N-4q_N\rint$  
\begin{equation}\label{mixtime4qlength}
d_{[x+1,x+4q_N],k}^{\go}(T)\le N^{-3}.
\end{equation}
\end{lemma}

 We are going to use the censoring inequality to guide all the particles to the right with the following plan.
We are going to design our censoring such that on the time interval $[iT, (i+1)T)$, with $i\in \bbZ_+$ satisfying $2(i+2 )q_N < N$, our  $k$ particles perform the exclusion process restricted in the interval on the interval  $\lint  2iq_N+1, 2(i+2 )q_N\rint$ (of length $4q_N$).
Hence at each such time step, particles take a time $T$ to shift towards the right of an amount $2q_N$.  After the whole $\lceil N/(2q_N)\rceil-1$ steps have been performed, all particles are in  $\lint N-4q_N+1, N\rint$.
Once this is done we conclude using censoring again by showing that the dynamics in $ \lint N-4q_N+1,N\rint$ with less than $q_N$ particles hits $\xi_{\max}$ after time $T$ with a positive probability.  For this last step we need the following result  which states roughly that $\xi_{\max}$  has a positive weight with high probability under the equilibrium measure.

 \begin{lemma} \label{lema:maxconfmassprim}
We have
\begin{equation}
\lim_{\gep \to 0} \inftwo{N\ge 1}{k\in \lint 1,N/2\rint} \bbP \left[ \pi^{\go}_{N,k}(\xi_{\max})>\gep \right]
=1.
\end{equation}
In particular if 
$\cB_{N,k} \colonequals \left\{\go: \pi_{[N-4q_N+1,N],k}^{\go}(\xi_{\max}) \ge 2 q_N^{-1} \right \},$
 we have
 \begin{equation}\label{event:cBN}
 \lim_{N \to \infty} \inf_{k \in \lint 1, q_N\rint} \bbP \left[ \cB_{N,k}\right]=1.
 \end{equation}
\end{lemma}

\begin{proposition}\label{prop:smallk}
If $k\le q_N$, if $\go \in \cA_N\cap \cB_{N,k}$ and setting
 $t_0:= T \left(\left\lceil \frac{N}{2q_N}  \right\rceil-1 \right)$, we have
\begin{equation}
  P_{t_0}(\xi_{\min},\xi_{\max}) \ge \frac{ 3}{2q_N}.
\end{equation}
 In particular   the inequality holds with high probability.
 
\end{proposition}
 The last part of the statement is of course a direct consequence of the first part combined with  \eqref{event:cBN} and of Proposition \ref{lema:createtrap} (which ensures that $\cA_N$ and $\cB_{N,k}$ are high probability events).
Before proving a proof of Proposition \ref{prop:smallk} using the strategy exposed above, let us use it to conclude the proof of the upper bound on the mixing time.
\begin{proof}[Proof of Theorem \ref{th:mixingub} when $k\le q_N$.] 
By Proposition \ref{prop:reduhit} and Proposition \ref{prop:smallk}, we have
\begin{equation}
d_{N,k}^{\go}(2q_Nt_0) \le \left(1- P_{t_0}(\xi_{\min},\xi_{\max}) \right)^{2q_N} \le \left(1-\frac{1}{q_N}\right)^{2q_N}\le \frac{1}{4},
\end{equation}
which allows us to conclude the proof for the case $k \le q_N$ with  the inequality
$$ \binom{4q_N}{q_N} \le \left(\frac{4^4}{3^3}\right)^{q_N}. $$

\end{proof}

\noindent Now we move to prove Proposition \ref{prop:smallk}  using the censoring inequality (Proposition \ref{PWcoro}).
 More precisely, we define for $i \in \lint  0, \lceil N/(2q_N) \rceil -3 \rint$
 \begin{equation}
  \cC_i \colonequals  \Big\{\{i2q_N,i2q_N+1\},    \{(i+2)2q_N,(i+2)2q_N+1\} \Big\}  
 \end{equation}
and set 
\begin{equation}
 \cC_{\lceil N/(2q_N) \rceil-2} \colonequals \left\{ N-4q_N, N-4q_N+1 \right\}. 
\end{equation}
We define a censoring scheme by setting 
\begin{equation}
\cC(t):= \cC_i   \quad \text{ for } t\in [i T,(i+1)T),  i\in \lint 0, \lceil N/(2q_N) \rceil-2\rint,
\end{equation}
and $\cC(t)=\emptyset$ for $t\ge \lceil N/(2q_N) \rceil-1$.
Let us write 
\begin{equation}
A_{\mathrm{ fin}}\colonequals \left\{ \xi \in \gO_{N,k} \ : \ \forall x\in \lint 0 ,N-4q_N \rint, \xi(x)=0  \right\}.
\end{equation}
Recalling the notation of Section \ref{compcensor}, we let $(\sigma^{{\min},\cC}_t)_{t \ge 0}$ denote the corresponding censored dynamics with initial condition $\xi_{\min}$. 
\begin{lemma}\label{lema:secondlast}
If $\go \in \cA_N$, we have 
\begin{equation}
 \bP \left[ \sigma^{{\min},\cC}_{(\lceil N/2q_N \rceil-2)T} \in A_{\mathrm{ fin}}\right] \ge  1-N^{-1}.
\end{equation}

\end{lemma}

\begin{proof}
For $i \in  \lint 0, \lceil N/2q_N \rceil -2\rint $,  we define 
\begin{equation*}
A_i \colonequals
 \left\{\xi \in \gO_{N,k}: 2iq_N< \bar \xi(1)  \le  \bar \xi(k) \le 2(i+2)q_N   \right\}.
\end{equation*}
Now we prove by induction that for all $i \in  \lint 0, \lceil N/2q_N \rceil -3 \rint$
\begin{equation}\label{inductmove}
\begin{split}
\bP \left[   \sigma^{{\min},\cC}_{iT} \in A_{i}   \right] \ge 1-i\frac{4q_N^2}{N^3}.
\end{split}
\end{equation}
From the definitions of $\cC$ and $\xi_{\min}$, the inequality in \eqref{inductmove} holds for $i=0$. Assuming that \eqref{inductmove} holds for $i$,  then $k$ particles perform the simple exclusion process restricted in the interval $\lint 2iq_N+1, 2(i+2)q_N \rint$. 
By  Lemma \ref{lema:LhalfRone} and Lemma \ref{lema:mixN3} with $x=2iq_N$, we have 
\begin{multline}\label{inductend}
\bP \left[ \sigma^{{\min},\cC}_{(i+1)T} \in A_{i+1}    \right]\\
 \ge \bP \left[   \sigma^{{\min},\cC}_{iT} \in A_i    \right]
- \left(  
\pi_{[2iq_N+1, 2(i+2)q_N],k}^{\go}\left(\bar \xi(1) \le 2(i+1)q_N \right) + d_{[2iq_N+1, 2(i+2)q_N],k}^{\go}(T)  \right)\\
\ge 1-i\frac{4q_N^2}{N^3}- \frac{4q_N^2}{N^3}.
\end{multline}
 This concludes the induction and the case 
  $i= \lceil N/2q_N \rceil -3 $ in \eqref{inductend} 
  to concludes the proof of the lemma since $$2(\lceil N/2q_N \rceil -2)q_N \ge N-4q_N.$$

\end{proof}

\begin{proof}[Proof of Proposition \ref{prop:smallk}]
 Using Proposition \ref{PWcoro}, it is sufficient to bound the corresponding probability for the censored dynamics, that is, $P^{\cC}_{  (\lceil N/2q_N \rceil-1)T}(\xi_{\min},\xi_{\max})$.
If $\sigma^{{\min},\cC}_{(\lceil N/2q_N \rceil-2)T} \in A_{\mathrm{fin}}$, then the restriction to the segment
$\lint N-4q_N+1,N\rint$ of the dynamics corresponds to an exclusion process with $k$ particles on a segment of length $4q_N$.
Let $\pi_{[N-4q_N+1,N],k}$  and $d_{[N-4q_N+1,N],k}(t)$ denote respectively the equilibrium measure and the distance to equilibrium for this dynamics, and then we have
\begin{multline}\label{lwbd:hitprobcensor}
P^{\cC}_{  (\lceil N/2q_N \rceil-1)T}(\xi_{\min},\xi_{\max}) \ge  \bP[ \sigma^{\xi_{\min},\cC}_{(\lceil N/2q_N \rceil-2)T} \in A_{\mathrm{fin}}](\pi_{[N-4q_N+1,N],k}(\xi_{\max})- d_{[N-4q_N+1,N],k}(T))\\
 \ge (1-N^{-1})(2q_N^{-1}-N^{-3}) \ge \frac{3}{2q_N}
\end{multline}
where we have used  the definition of $\cB_{N,k}$ (recall \eqref{event:cBN}) and Lemma \ref{lema:mixN3} with $x=N-4q_N$.

\end{proof}

\subsection{Proof of auxiliary lemmas}

\begin{proof}[Proof of Lemma \ref{lema:LhalfRone}.] To provide an upper bound on $\pi_{[x+1,x+4q_N],k}^{\go} \left[ \bar \xi(1) \le x+2q_N \right]$, for $\xi \in \gO_{[x+1,x+4q_N],k}$ we define its rightmost empty site to be
\begin{equation}
\bar R(\xi) \colonequals \sup \left\{y \in \lint x+1,x+4q_N \rint: \ \xi(y)=0 \right\}.
\end{equation}
   As in \eqref{exchange01}, we have
\begin{equation}
\pi_{[x+1,x+4q_N],k}^{\go} \left[\bar \xi(1)=z,  \bar R(\xi)=y \right] \le e^{V^{\go}(y)-V^{\go}(z)} \le N^{-3}
\end{equation}
where we have used $y-z \ge q_N$ and $\go \in \cA_N$.
Then we have
\begin{multline}
\pi_{[x+1,x+4q_N],k}^{\go} \left[ \bar \xi(1) \le x+2q_N\right]= \sumtwo{z \in \lint x+ 1, x+2q_N\rint}{y \in \lint x+4q_N-k+2,x+4q_N\rint}\pi_{[x+1,x+4q_N],k} \left[\bar \xi(1)=z,  \bar R(\xi)=y \right]\\
\le 2q_N^2 N^{-3}.
\end{multline}

We now move to deal with $\pi_{[x+1,x+4q_N],q_N}^{\go} \left[ \xi(x+4q_N)=0 \right]$. For $\xi \in \gO_{[x+1,x+4q_N],q_N}$, we define its leftmost particle to be
$$ \bar L(\xi) \colonequals \inf \left\{y \in \lint x+1, x+4q_N\rint: \ \xi(y)=1 \right\}. $$
As in \eqref{exchange01}, we have
\begin{equation}
\pi_{[x+1,x+4q_N],q_N}^{\go} \left[ \xi(x+4q_N)=0; \bar L(\xi)=y \right] \le e^{V^{\go}(x+4q_N)-V^{\go}(y)} \le N^{-3}
\end{equation}
where we have used $y \le x+3q_N$ and $\go \in \cA_N$. Then
\begin{multline}
\pi_{[x+1,x+4q_N],q_N}^{\go} \left[ \xi(x+4q_N)=0\right]= \sum_{y \in \lint x+1, x+3q_N \rint}\pi_{[x+1,x+4q_N],q_N}^{\go} \left[ \xi(x+4q_N)=0; \bar L(\xi)=y \right] \\
\le 3q_NN^{-3}.
\end{multline}

\end{proof}

\begin{proof}[Proof of Lemma \ref{lema:maxconfmassprim}]
Recall the event $\cA_r$ in \eqref{event:Ar}. Observe that
\begin{equation}
\max_{\xi \in \cA_r} \left( V^{\go}(\xi_{\max})-V^{\go}(\xi) \right) \le 2r^2 \log \frac{1-\ga}{\ga},
\end{equation}
and then we have
\begin{equation}\label{maxconfratio}
\frac{\pi_{N,k}^{\go}(\xi_{\max})}{\pi_{N,k}^{\go}(\cA_r)} \ge \vert \cA_r \vert^{-1} \exp\left(-\max_{\xi \in \cA_r} \left( V^{\go}(\xi_{\max})-V^{\go}(\xi) \right) \right) \ge 2^{-2r} e^{-2r^2 \log \frac{1-\ga}{\ga}}.
\end{equation}
For given $\gep>0$ sufficiently small, we take 
\begin{equation}\label{chooser}
r(\gep) \colonequals \left(\frac{-\log 2\gep }{ \log \frac{2(1-\ga)}{\ga}} \right)^{1/2} 
\end{equation}
so that the rightmost hand-side of \eqref{maxconfratio} is larger than $2\gep$.
Moreover, by \eqref{prob:rempty} we know that
\begin{equation}\label{prob:eventAr} 
  \lim_{r\to \infty} \inftwo{N\ge 1}{ k\in \lint 1, N/2\rint}\bbP\left[ \pi_{N,k}^{\go} \left( \cA_r \right) \ge 1-2 (1-e^{\frac { \bbE[\log \rho_1]} 2})^{-2}e^{\frac{{ \bbE[\log \rho_1]} r}2} \right]=1.
\end{equation}
Since  when $r$ is sufficiently large we have $$1-2 (1-e^{\frac { \bbE[\log \rho_1]} 2})^{-2}e^{\frac{{ \bbE[\log \rho_1]} r}2} \ge \frac{1}{2},$$  then by \eqref{prob:eventAr} with $r$ chosen as in \eqref{chooser} we obtain
\begin{equation}
\lim_{\gep \to 0} \inftwo{N\ge 1}{ k\in \lint 1, N/2\rint}\bbP\left[ 
\pi_{N,k}^{\go}(\xi_{\max}) \ge \gep \right]=1.
\end{equation}

\end{proof}
\subsection{The case $k_N\ge q_N$}
\label{morethank}

To treat the case of a larger number of particles, the small problem there is with the strategy of the previous subsection is that it does not allow to channel all the $k$ particles to the right at the same time. What we do instead  is that we use the process to transport one particle to the right, and then use Proposition \ref{prop:censor} to be able  to move all other particles to the left and iterate the process.
We largely recycle the strategy used in the previous section.  
 In the final step as in \eqref{lwbd:hitprobcensor}, we need to deal with
 the leftmost $q_N$ particles performing the exclusion process restricted in the interval $\lint N-k-3q_N+1,  N-k+q_N\rint $, and then define
  $$\cB'_{N,k}=\left\{\go: \pi_{[N-k-3q_N+1,N-k+q_N],q_N}^{\go}(\xi_{\max}') \ge 2q_N^{-1} \right \}$$
where $\xi_{\max}' \colonequals \ind_{\{N-k+1 \le  x \le N-k+q_N\}}$. By Lemma \ref{lema:maxconfmassprim} 
we have
\begin{equation}
 \lim_{N \to \infty} \inf_{k \in \lint q_N+1,N/2\rint}\bbP \left[ \cB_{N,k}' \right]=1.
\end{equation}

\begin{proposition}\label{lastproposition}
 If $k> q_N$ and $\go \in \cA_N\cap \cB'_{N,k}$,  setting
 $t_1 \colonequals \left( \left \lceil  \frac{N-k+q_N}{2q_N} \right \rceil-1 \right)(k-q_N+1 )T $ we have
\begin{equation}
  P_{t_1}(\xi_{\min},\xi_{\max}) \ge \frac{1}{q_N}.
\end{equation}

\end{proposition}

\begin{proof}[Proof of Theorem \ref{th:mixingub} when $k> q_N$]
By Proposition \ref{prop:reduhit} and Proposition \ref{lastproposition}, we have
\begin{equation}
d_{N,k}^{\go}(2q_Nt_1) \le \left(1- P_{t_1}(\xi_{\min},\xi_{\max}) \right)^{2q_N} \le \left(1-\frac{1}{q_N}\right)^{2q_N}\le \frac{1}{4},
\end{equation}
which allows us to conclude the proof for the case $k > q_N$ with  the inequality
$$ \binom{4q_N}{q_N} \le \left(\frac{4^4}{3^3}\right)^{q_N}. $$

\end{proof}

 The remaining of the subsection is devoted to the proof of Proposition \ref{lastproposition}.
This time we need to combine our censoring scheme with displacements of particles to the  left (using  Proposition \ref{prop:censor}).
Our plan is to first move (one by one) the rightmost $k-q_N$ particles  to the segment $\lint N-k+q_N+1,  N\rint$ and use censoring to block the these $k-q_N$ particles  afterwards. 
We are then left with the problem of moving the remaining $q_N$ particles, and this can be treated as in Proposition \ref{prop:smallk}.

\medskip

Let us explain our plan to move the the rightmost $k-q_N$ particles one by one with censoring and displacement. We proceed by induction (each step is going to leave aside an event of small probability, and our technical estimates are such that the sum over all steps of these probabilities will remain small).
We set $r= \lceil (N-k+q_N)/2q_N \rceil -1$,
 and define for $j\in \lint 0, k-q_N\rint$, $i\in \lint  0, \lceil (N-k+q_N)/2q_N \rceil -3 \rint$,  $a_{i,j} \colonequals k-q_N-j+2q_Ni$,
 \begin{equation}\begin{split}\label{choicefocensoring}
  \cC_{i,j}&:=  \Big\{\{a_{i,j},a_{i,j}+1\},    
  \{a_{i,j}+4q_N,a_{i,j}+4q_N+1\}, \{N-j,N-j+1\}\Big\}  ,\\
   \cC^*_{j}&=  \Big\{\{ N-4q_N-j,N-4q_N-j+1\},    
  \{ N-j, N-j+1\} \Big\}.
   \end{split}
 \end{equation}
We define the censoring  scheme $\cC$ by setting
\begin{equation}\label{secondcensor}
 \begin{cases}
  \cC(t)=   \cC_{i,j}  &\text{ if }  t\in[(i+rj)T, (i+rj+1)T ),\\
  \cC(t)=   \cC^*_{j}  &\text{ if }  t\in[\left(r(j+1)-1\right) T, r(j+1)T ),\\
\cC(t)=\emptyset &\text{ if } t\ge r(k-q_N+1)T.
  \end{cases}
\end{equation}

The censored dynamic $(\sigma^{\cC,\min}_t)$ moves the first particle to the right in a time $rT$. Indeed, the same mechanism used in the proof of Proposition \ref{prop:smallk} moves (w.h.p) the last $q_N$ particles in the segment $\lint N-4q_N+1, N\rint$ by time $(r-1)T$. Then we mix the $q_N$ particles within the  segment $\lint N-4q_N+1, N\rint$ and Lemma \ref{lema:LhalfRone} ensures that after an additional time $T$, the last site $N$ is occupied by a particle.

\medskip

We then proceed by induction to show that for $j\le k-q_N$ all the sites in the segment $\lint N-j+1,  N \rint$ are occupied by particles by time $rjT$. Our censoring is designed so that 
after time $rjT$ the number of particles in the $j$ rightmost sites does not change.

\medskip

In order to facilitate the induction (this is not strictly necessary though) at each time of the form $rjT=:s_j$ we move all the leftmost $N-j$ particles to the left on the segment $\lint 1, N-j\rint$, so that the beginning of each induction step looks the same.
We define thus $Q_j$   by setting 
\begin{equation}
 Q_j(\xi, \xi_j^*)=1, \quad Q_j(\xi, \xi')=0 \text{ if } \xi'\ne \xi^*_j
\end{equation}
where the function $\xi\to \xi^*_j$ is defined by (recall \eqref{partiposition})
\begin{equation} \label{displaceconfig}
\bar \xi^*_{j}(\ell)= \begin{cases}   \ell &\quad \text{ if } l \le  k-j,\\
                   \bar \xi(\ell) &\quad \text{ if } \ell>k-j.
                  \end{cases}
    \end{equation}    
Since $\xi^*_j\le \xi$, $Q_j$ satisfies \eqref{upp}.
We  let $(\tilde \sigma_t)_{t \ge 0}$ denote the composed censored dynamics (recall \eqref{semigroup:displace}) corresponding to $\cC$, $(s_j)_{j=1}^{k-q_N}$ and $(Q_j)_{j=1}^{ k-q_N}$  and starting from $\xi_{\min}$.
We set 
$$\xi^0_{j}:= \ind_{\lint 1, k-j\rint}+ \ind_{\lint N-j+1, N\rint}.$$
The following lemma formalizes in a quantitative manner the induction described above.

\medskip

\begin{lemma}\label{lema:inductmoveone}
 For all  $j \in \lint 0, k-q_N\rint$, we have
 \begin{equation}
  \bP\left[  \tilde \sigma_{rjT}=\xi^0_j \right] \ge 1- 4jq_N N^{-2}.
\end{equation}
 
\end{lemma}

\begin{proof}
 The statement is trivial for $j=0$.
 For the induction step it is sufficient to prove that 
 \begin{equation}
    \bP\left[  \tilde \sigma_{r(j+1)T}=\xi^0_{j+1} \ | \   \tilde \sigma_{rjT}=\xi^0_{j}\right] \ge 1- 4q_N N^{-2}.
 \end{equation}
With our choice for $\cC$, the $j$ particles in the interval $\lint N-j+1, N\rint$ do not move  
between time instants $rjT$ and $r(j+1)T$, it is therefore sufficient to control 
  $\bP\left[  \tilde \sigma_{r(j+1)T}(N-j)=1 \ | \   \tilde \sigma_{rjT}=\xi^0_{j}\right]$.
 Let us define
\begin{equation}
B_j \colonequals \left\{\xi \in \gO_{N,k}: \sum_{N-j-4q_N+1}^{N-j}\xi(x)=q_N \right\}
\end{equation}
We can repeat the proof of Lemma \ref{lema:secondlast} to obtain that
\begin{equation}\label{repeatinduct}
\bP \left[ \tilde \gs_{rjT+(r-1)T}  \in B_j \ \vert  \ \tilde \sigma_{rjT}=\xi^*_j  \right] \ge 1-(r-1)\frac{4q_N^2}{N^3}.
\end{equation}
Now in the time interval $[rjT+(r-1)T, r(j+1)T)$, the censoring makes the restriction of the dynamics to the segment  $\lint N-j-4q_N+1, N-j\rint$ an exclusion process with $q_N$ particles.
Hence using  Lemma \ref{lema:mixN3}  and the second estimate in Lemma \ref{lema:LhalfRone} we have
for any $\chi\in B_j$
\begin{equation}\label{packonepart}
 \bP \left[ \tilde \gs_{r(j+1)T}(N-j)=1 \ \vert  \ \tilde \gs_{rjT+(r-1)T}  = \chi \right]
 \ge 1- N^{-3}(1+3q^2_N).
\end{equation}
Combining \eqref{repeatinduct} and \eqref{packonepart}, we obtain
\begin{equation}
\bP\left[   \tilde \sigma_{r(j+1)T}=\xi^0_{j}\right] \ge \bP\left[   \tilde \sigma_{r(j+1)T}=\xi^0_{j}\right]-r\frac{4q_N^2}{N^3} \ge 1-4(j+1)q_N N^{-2}.
\end{equation}
\end{proof}

\begin{proof}[Proof of Proposition \ref{lastproposition}]

  Taking $j=k-q_N$ in Lemma \ref{lema:inductmoveone}, from now on we assume that the event $\{ \tilde\gs_{(k-q_N)rT}= \xi_{k-q_N}^0\}$ holds. Then the rightmost $k-q_N$ particles are frozen in the rightmost $k-q_N$ sites for $t \ge (k-q_N)rT $, and 
at $t = (k-q_N)rT $ the leftmost $q_N$ particles are in the leftmost $q_N$ sites. Thus we can repeat the proof in Proposition \ref{prop:smallk} to obtain  
  \begin{equation}
\bP\left[  \tilde \sigma_{r(k-q_N+1)T}=\xi_{\max} \right] \ge \frac{3}{2} q_N^{-1} \left(1- (k-q_N) \frac{4q_N}{N^2} \right)   \ge \frac{1}{q_N}
\end{equation}
where we have used $\go \in \cB_{N,k}'$.
We conclude the proof by Proposition \ref{PWcoro} and Proposition \ref{prop:censor}.
\end{proof}

\bibliographystyle{alpha}
\bibliography{library.bib}

\begin{thebibliography}{BBHM05}

\bibitem[AFJV15]{avena2015symmetric}
Luca Avena, Tertuliano Franco, Milton Jara, and Florian V\"{o}llering.
\newblock Symmetric exclusion as a random environment: hydrodynamic limits.
\newblock {\em Ann. Inst. Henri Poincar\'{e} Probab. Stat.}, 51(3):901--916,
  2015.

\bibitem[Ald83]{Aldous83}
David Aldous.
\newblock Random walks on finite groups and rapidly mixing {M}arkov chains.
\newblock In {\em Seminar on probability, {XVII}}, volume 986 of {\em Lecture
  Notes in Math.}, pages 243--297. Springer, Berlin, 1983.

\bibitem[BBHM05]{benjamini2005mixing}
Itai Benjamini, Noam Berger, Christopher Hoffman, and Elchanan Mossel.
\newblock Mixing times of the biased card shuffling and the asymmetric
  exclusion process.
\newblock {\em Transactions of the American Mathematical Society},
  357(8):3013--3029, 2005.

\bibitem[BECE00]{blythe2000exact}
RA~Blythe, MR~Evans, F~Colaiori, and FHL Essler.
\newblock Exact solution of a partially asymmetric exclusion model using a
  deformed oscillator algebra.
\newblock {\em Journal of Physics A: Mathematical and General}, 33(12):2313,
  2000.

\bibitem[DSC93]{DSC93}
Persi Diaconis and Laurent Saloff-Coste.
\newblock Comparison theorems for reversible {M}arkov chains.
\newblock {\em Ann. Appl. Probab.}, 3(3):696--730, 1993.

\bibitem[Fag08]{faggionato2008random}
Alessandra Faggionato.
\newblock Random walks and exclusion processes among random conductances on
  random infinite clusters: homogenization and hydrodynamic limit.
\newblock {\em Electronic Journal of Probability}, 13:2217--2247, 2008.

\bibitem[{Fag}20]{Faggionato2020}
Alessandra {Faggionato}.
\newblock {Hydrodynamic limit of simple exclusion processes in symmetric random
  environments via duality and homogenization}.
\newblock {\em arXiv e-prints}, page arXiv:2011.11361, November 2020.

\bibitem[FGS16]{FGS16}
Tertuliano Franco, Patr\'{\i}cia Gon\c{c}alves, and Marielle Simon.
\newblock Crossover to the stochastic {B}urgers equation for the {WASEP} with a
  slow bond.
\newblock {\em Comm. Math. Phys.}, 346(3):801--838, 2016.

\bibitem[FN17]{FN17}
Tertuliano Franco and Adriana Neumann.
\newblock Large deviations for the exclusion process with a slow bond.
\newblock {\em Ann. Appl. Probab.}, 27(6):3547--3587, 2017.

\bibitem[FRS19]{FRS19}
Simone {Floreani}, Frank {Redig}, and Federico {Sau}.
\newblock {Hydrodynamics for the partial exclusion process in random
  environment}.
\newblock {\em arXiv e-prints}, page arXiv:1911.12564, November 2019.

\bibitem[GK13]{gantert2012cutoff}
Nina Gantert and Thomas Kochler.
\newblock Cutoff and mixing time for transient random walks in random
  environments.
\newblock {\em ALEA Lat. Am. J. Probab. Math. Stat.}, 10(1):449–484, 2013.

\bibitem[GNS20]{gantertboundaries}
Nina {Gantert}, Evita {Nestoridi}, and Dominik {Schmid}.
\newblock {Mixing times for the simple exclusion process with open boundaries}.
\newblock {\em arXiv e-prints}, page arXiv:2003.03781, March 2020.

\bibitem[HKT20]{hilario2020random}
Marcelo~R Hil{\'a}rio, Daniel Kious, and Augusto Teixeira.
\newblock Random walk on the simple symmetric exclusion process.
\newblock {\em Communications in Mathematical Physics}, 379(1):61--101, 2020.

\bibitem[HS15]{huveneers2015random}
Fran{\c{c}}ois Huveneers and Fran{\c{c}}ois Simenhaus.
\newblock Random walk driven by simple exclusion process.
\newblock {\em Electronic Journal of Probability}, 20, 2015.

\bibitem[Jar11]{jara2011hydrodynamic}
Milton Jara.
\newblock Hydrodynamic limit of the exclusion process in inhomogeneous media.
\newblock In {\em Dynamics, Games and Science II}, pages 449--465. Springer,
  2011.

\bibitem[JM20]{jara2020symmetric}
Milton Jara and Ot{\'a}vio Menezes.
\newblock Symmetric exclusion as a random environment: invariance principle.
\newblock {\em Annals of Probability}, 48(6):3124--3149, 2020.

\bibitem[Kin73]{kingman73}
J.~F.~C. Kingman.
\newblock Subadditive ergodic theory.
\newblock {\em Ann. Probability}, 1:883--909, 1973.

\bibitem[KKS75]{kesten1975limit}
Harry Kesten, Mykyta~V Kozlov, and Frank Spitzer.
\newblock A limit law for random walk in a random environment.
\newblock {\em Compositio Mathematica}, 30(2):145--168, 1975.

\bibitem[KOV89]{KOV89}
C.~Kipnis, S.~Olla, and S.~R.~S. Varadhan.
\newblock Hydrodynamics and large deviation for simple exclusion processes.
\newblock {\em Comm. Pure Appl. Math.}, 42(2):115--137, 1989.

\bibitem[Lac16a]{lacoinprofile}
Hubert Lacoin.
\newblock The cutoff profile for the simple exclusion process on the circle.
\newblock {\em Ann. Probab.}, 44(5):3399--3430, 2016.

\bibitem[Lac16b]{lacoin2016mixing}
Hubert Lacoin.
\newblock Mixing time and cutoff for the adjacent transposition shuffle and the
  simple exclusion.
\newblock {\em The Annals of Probability}, 44(2):1426--1487, 2016.

\bibitem[Lig12]{liggett2012interacting}
Thomas~Milton Liggett.
\newblock {\em Interacting particle systems}, volume 276.
\newblock Springer Science \& Business Media, 2012.

\bibitem[LL19]{labbe2016cutoff}
Cyril Labb{\'e} and Hubert Lacoin.
\newblock Cutoff phenomenon for the asymmetric simple exclusion process and the
  biased card shuffling.
\newblock {\em The Annals of Probability}, 47(3):1541--1586, 2019.

\bibitem[LL20]{labbe2018mixing}
Cyril Labb{\'e} and Hubert Lacoin.
\newblock Mixing time and cutoff for the weakly asymmetric simple exclusion
  process.
\newblock {\em Annals of Applied Probability}, 30(4):1847--1883, 2020.

\bibitem[LP16]{levin2016mixing}
David~A Levin and Yuval Peres.
\newblock Mixing of the exclusion process with small bias.
\newblock {\em Journal of Statistical Physics}, 165(6):1036--1050, 2016.

\bibitem[LP17]{LPWMCMT}
David~A. Levin and Yuval Peres.
\newblock {\em Markov chains and mixing times}.
\newblock American Mathematical Society, Providence, RI, 2017.
\newblock Second edition of [ MR2466937], With contributions by Elizabeth L.
  Wilmer, With a chapter on ``Coupling from the past'' by James G. Propp and
  David B. Wilson.

\bibitem[Mor06]{Morris06}
Ben Morris.
\newblock The mixing time for simple exclusion.
\newblock {\em Ann. Appl. Probab.}, 16(2):615--635, 2006.

\bibitem[PW13]{peres2013can}
Yuval Peres and Peter Winkler.
\newblock Can extra updates delay mixing?
\newblock {\em Communications in Mathematical Physics}, 323(3):1007--1016,
  2013.

\bibitem[Qua92]{Quastel92}
Jeremy Quastel.
\newblock Diffusion of color in the simple exclusion process.
\newblock {\em Comm. Pure Appl. Math.}, 45(6):623--679, 1992.

\bibitem[Rez91]{rez91}
Fraydoun Rezakhanlou.
\newblock Hydrodynamic limit for attractive particle systems on {${\bf Z}^d$}.
\newblock {\em Comm. Math. Phys.}, 140(3):417--448, 1991.

\bibitem[Ros81]{Rost81}
H.~Rost.
\newblock Nonequilibrium behaviour of a many particle process: density profile
  and local equilibria.
\newblock {\em Z. Wahrsch. Verw. Gebiete}, 58(1):41--53, 1981.

\bibitem[Sch19]{schmid2019mixing}
Dominik Schmid.
\newblock Mixing times for the simple exclusion process in ballistic random
  environment.
\newblock {\em Electronic Journal of Probability}, 24, 2019.

\bibitem[Sin82]{sinai1982limit}
Yakov~Grigor'evich Sinai.
\newblock Limit behaviour of one-dimensional random walks in random
  environments.
\newblock {\em Teoriya Veroyatnostei i ee Primeneniya}, 27(2):247--258, 1982.

\bibitem[Sol75]{solomon1975random}
Fred Solomon.
\newblock Random walks in a random environment.
\newblock {\em The annals of probability}, pages 1--31, 1975.

\bibitem[Szn04]{sznitman2004topics}
Alain-Sol Sznitman.
\newblock Topics in random walks in random environment.
\newblock In {\em School and conference on probability theory: 13-17 May 2002},
  volume~17, pages 203--266. The Abdus Salam International Centre for
  Theoretical Physics, 2004.

\bibitem[Wil04]{wilson2004mixing}
David~Bruce Wilson.
\newblock Mixing times of lozenge tiling and card shuffling markov chains.
\newblock {\em The Annals of Applied Probability}, 14(1):274--325, 2004.

\bibitem[Yau97]{Yau97}
Horng-Tzer Yau.
\newblock Logarithmic {S}obolev inequality for generalized simple exclusion
  processes.
\newblock {\em Probab. Theory Related Fields}, 109(4):507--538, 1997.

\bibitem[Zei04]{zeitouni2004part}
Ofer Zeitouni.
\newblock Part ii: Random walks in random environment.
\newblock In {\em Lectures on Probability Theory and Statistics}, pages
  190--312. Springer, 2004.

\end{thebibliography}

\end{document}